\documentclass[12pt]{article}
\usepackage{amsmath,amssymb,amsfonts,amsthm}
\usepackage{eucal} 
\usepackage{pstricks,pst-plot,pst-node}
\usepackage[dvips]{graphics}

\usepackage{texdraw}
\input{txdtools}

\def\fatness{1.2}
\everytexdraw{%
	\drawdim pt \linewd 0.5 \textref h:C v:C
  \arrowheadsize l:4 w:3  
  \arrowheadtype t:V
}

\setbox7=\hbox{
\raisebox{1pt}{\begin{texdraw}
	\drawdim pt \linewd 0.3 \setunitscale 1.2
	\move (-3 0)
	\move (0 -1.2) 
	\rlvec (6 0) 
	\rlvec (0 -1.2)
	\rlvec (2.4 2.4) 
	\rlvec (-2.4 2.4)
	\rlvec (0 -1.2) 
	\rlvec (-6 0)
	\rlvec (0 -2.4)
	\move (11.5 0)
\end{texdraw}}}

\newcommand{\idFA}{%
	\bsegment
		\realmult {\fatness} {2.5} \xrad
		\realmult {\fatness} {6.25} \yrad
		\move (0 0) \lellip rx:{\xrad} ry:{\yrad}
		\move (40 0) \lellip rx:{\xrad} ry:{\yrad}
		\realadd { 0} { \yrad} \inB
		\realadd { 0} {-\yrad} \inA
		\move (0 {\inB}) \rlvec (40 0) 
		\move (0 {\inA}) \rlvec (40 0) 
		\savepos (40 0)(*ex *ey)
	\esegment
	\move (*ex *ey)
}
\newcommand{\counitFA}{%
	\bsegment
		\realmult {\fatness} {2.5} \xrad
		\realmult {\fatness} {6.25} \yrad
		\move (0 0)
		\move (20 0)
		\move (0 0) \lellip rx:{\xrad} ry:{\yrad}
		\realadd { 0} { \yrad} \inB
		\realadd { 0} {-\yrad} \inA
		\move (0 {\inB}) \rlvec (8 0) 
		\move (0 {\inA}) \rlvec (8 0) 
	  \move (8 {\inA}) \clvec (20 {\inA})(20 {\inB})(8 {\inB}) 

		\savepos (0 0)(*ex *ey)
	\esegment
	\move (*ex *ey)
}
\newcommand{\unitFA}{%
	\bsegment
		\realmult {\fatness} {2.5} \xrad
		\realmult {\fatness} {6.25} \yrad
		\move (0 0)
		\move (-20 0)
		\move (0 0) \lellip rx:{\xrad} ry:{\yrad}
		\realadd { 0} { \yrad} \inB
		\realadd { 0} {-\yrad} \inA
		\move (0 {\inB}) \rlvec (-8 0) 
		\move (0 {\inA}) \rlvec (-8 0) 
	  \move (-8 {\inA}) \clvec (-20 {\inA})(-20 {\inB})(-8 {\inB}) 

		\savepos (0 0)(*ex *ey)
	\esegment
	\move (*ex *ey)
}
\newcommand{\twistFA}{%
	\bsegment
		\realmult {\fatness} {2.5} \xrad
		\realmult {\fatness} {6.25} \yrad
		\move (0 0) \lellip rx:{\xrad} ry:{\yrad}
		\move (0 30) \lellip rx:{\xrad} ry:{\yrad}
		\move (40 0) \lellip rx:{\xrad} ry:{\yrad}
		\move (40 30) \lellip rx:{\xrad} ry:{\yrad}
		\realadd {30} { \yrad} \inD
		\realadd {30} {-\yrad} \inC
		\realadd { 0} { \yrad} \inB
		\realadd { 0} {-\yrad} \inA
		\move (0 {\inD}) \clvec (28 {\inD})(20 {\inB})(40 {\inB}) 
		\move (0 {\inC}) \clvec (20 {\inC})(12 {\inA})(40 {\inA}) 

		\move (0 {\inB}) \clvec (20 {\inB})(12 {\inD})(40 {\inD}) 
		\move (0 {\inA}) \clvec (28 {\inA})(20 {\inC})(40 {\inC}) 
		\savepos (40 0)(*ex *ey)
	\esegment
	\move (*ex *ey)
}

\newcommand{\multFA}{%
	\bsegment
		\realmult {\fatness} {2.5} \xrad
		\realmult {\fatness} {6.25} \yrad
		\move (0 0) \lellip rx:{\xrad} ry:{\yrad}
		\move (0 30) \lellip rx:{\xrad} ry:{\yrad}
		\move (40 15) \lellip rx:{\xrad} ry:{\yrad}
		\realadd {15} { \yrad} \outB
		\realadd {15} {-\yrad} \outA
		\realadd {30} { \yrad} \inD
		\realadd {30} {-\yrad} \inC
		\realadd { 0} { \yrad} \inB
		\realadd { 0} {-\yrad} \inA
		\move (0 {\inD}) \clvec (25 {\inD})(20 {\outB})(40 {\outB}) 
		\move (0 {\inB}) \clvec (18 {\inB})(18 {\inC})(0 {\inC}) 
		\move (0 {\inA}) \clvec (25 {\inA})(20 {\outA})(40 {\outA}) 
		\savepos (40 15)(*ex *ey)
	\esegment
	\move (*ex *ey)
}

\newcommand{\comultFA}{%
	\bsegment
		\realmult {\fatness} {2.5} \xrad
		\realmult {\fatness} {6.25} \yrad
		\move (0 0) \lellip rx:{\xrad} ry:{\yrad}
		\move (40 -15) \lellip rx:{\xrad} ry:{\yrad}
		\move (40 15) \lellip rx:{\xrad} ry:{\yrad}
		\realadd { 0} { \yrad} \outB
		\realadd { 0} {-\yrad} \outA
		\realadd {15} { \yrad} \inD
		\realadd {15} {-\yrad} \inC
		\realadd {-15} { \yrad} \inB
		\realadd {-15} {-\yrad} \inA
		\move (40 {\inD}) \clvec (15 {\inD})(20 {\outB})(0 {\outB}) 
		\move (40 {\inB}) \clvec (22 {\inB})(22 {\inC})(40 {\inC}) 
		\move (40 {\inA}) \clvec (15 {\inA})(20 {\outA})(0 {\outA}) 
		\savepos (40 -15)(*ex *ey)
	\esegment
	\move (*ex *ey)
}

\newcommand{\SFA}{
	\bsegment
		\move ( 0  0) \lellip rx:3 ry:7.5
		\move ( 0 30) \lellip rx:3 ry:7.5
		\move (40 45) \lellip rx:3 ry:7.5
		\move (40 15) \lellip rx:3 ry:7.5
		\move (40 -15) \lellip rx:3 ry:7.5
		\move (80  0) \lellip rx:3 ry:7.5
		\move (80 30) \lellip rx:3 ry:7.5

    \move (0 37.5) \clvec (20 37.5)(15 52.5)(40 52.5)
                   \clvec (65 52.5)(60 37.5)(80 37.5) 
    \move (0 -7.5) \clvec (20 -7.5)(15 -22.5)(40 -22.5)
                   \clvec (65 -22.5)(60 -7.5)(80 -7.5) 
		\move (80 22.5) \clvec (55 22.5)(60 37.5)(40 37.5)
                    \clvec (22 37.5)(22 22.5)(40 22.5)
                    \clvec (65 22.5)(60 7.5)(80 7.5) 
		\move (0 7.5) \clvec (25 7.5)(20 -7.5)(40 -7.5)
                    \clvec (58 -7.5)(58 7.5)(40 7.5)
                    \clvec (15 7.5)(20 22.5)(0 22.5) 
		\savepos (80 0)(*ex *ey)
	\esegment
  \move (*ex *ey)
}

\newcommand{\onedot}[1]{
  \bsegment
	\move (0 0) \fcir f:0 r:2
	\esegment
}

\newcommand{\boundaryFA}{%
  \bsegment
    \realmult {\fatness} {2.5} \xrad
    \realmult {\fatness} {6.25} \yrad
    \move (0 0) \lellip rx:{\xrad} ry:{\yrad}
    \savepos (0 0)(*ex *ey)
  \esegment
  \move (*ex *ey)
}

\newcommand{\multNB}{%
	\bsegment
		\realmult {\fatness} {2.5} \xrad
		\realmult {\fatness} {6.25} \yrad
		\move (0 0) 
		\move (0 30) 
		\move (40 15) 
		\realadd {15} { \yrad} \outB
		\realadd {15} {-\yrad} \outA
		\realadd {30} { \yrad} \inD
		\realadd {30} {-\yrad} \inC
		\realadd { 0} { \yrad} \inB
		\realadd { 0} {-\yrad} \inA
		\move (0 {\inD}) \clvec (25 {\inD})(20 {\outB})(40 {\outB}) 
		\move (0 {\inB}) \clvec (18 {\inB})(18 {\inC})(0 {\inC}) 
		\move (0 {\inA}) \clvec (25 {\inA})(20 {\outA})(40 {\outA}) 
		\savepos (40 15)(*ex *ey)
	\esegment
	\move (*ex *ey)
}

\newcommand{\comultNB}{%
	\bsegment
		\realmult {\fatness} {2.5} \xrad
		\realmult {\fatness} {6.25} \yrad
		\move (0 0) 
		\move (40 -15) 
		\move (40 15) 
		\realadd { 0} { \yrad} \outB
		\realadd { 0} {-\yrad} \outA
		\realadd {15} { \yrad} \inD
		\realadd {15} {-\yrad} \inC
		\realadd {-15} { \yrad} \inB
		\realadd {-15} {-\yrad} \inA
		\move (40 {\inD}) \clvec (15 {\inD})(20 {\outB})(0 {\outB}) 
		\move (40 {\inB}) \clvec (22 {\inB})(22 {\inC})(40 {\inC}) 
		\move (40 {\inA}) \clvec (15 {\inA})(20 {\outA})(0 {\outA}) 
		\savepos (40 -15)(*ex *ey)
	\esegment
	\move (*ex *ey)
}

\newdimen\tableauside\tableauside=1.0ex
\newdimen\tableaurule\tableaurule=0.4pt
\newdimen\tableaustep
\def\phantomhrule#1{\hbox{\vbox to0pt{\hrule height\tableaurule width#1\vss}}}
\def\phantomvrule#1{\vbox{\hbox to0pt{\vrule width\tableaurule height#1\hss}}}
\def\sqr{\vbox{%
  \phantomhrule\tableaustep
  \hbox{\phantomvrule\tableaustep\kern\tableaustep\phantomvrule\tableaustep}%
  \hbox{\vbox{\phantomhrule\tableauside}\kern-\tableaurule}}}
\def\squares#1{\hbox{\count0=#1\noindent\loop\sqr
  \advance\count0 by-1 \ifnum\count0>0\repeat}}
\def\tableau#1{\vcenter{\offinterlineskip
  \tableaustep=\tableauside\advance\tableaustep by-\tableaurule
  \kern\normallineskip\hbox
    {\kern\normallineskip\vbox
      {\gettableau#1 0 }%
     \kern\normallineskip\kern\tableaurule}%
  \kern\normallineskip\kern\tableaurule}}
\def\gettableau#1 {\ifnum#1=0\let\next=\null\else
  \squares{#1}\let\next=\gettableau\fi\next}

\newpsobject{showgrid}{psgrid}{subgriddiv=1,griddots=5,gridlabels=0pt}

\newcommand{\lang}{\left\langle}
\newcommand{\rang}{\right\rangle}

\newcommand{\brang}{\big\rangle}

\newcommand{\Bv}{\Big |}
\newcommand{\bv}{\big |}
\newcommand{\lv}{\left |}

\newcommand{\zz}{{\mathfrak{z}}}
\newcommand{\znums}{{\mathbb Z}}
\newcommand{\com}{{\mathbb C}}
\newcommand{\vac}{v_\emptyset}

\newcommand{\bS}{\mathsf{S}}

\newcommand{\bJ}{\mathsf{J}}
\newcommand{\TT}{T^{\pm}}

\newcommand{\cF}{\mathcal{F}}
\newcommand{\bW}{\mathsf{W}}
\newcommand{\bw}{\mathsf{w}}

\newcommand{\MM}{\mathsf{M}}

\newcommand{\fpa}{0_\uparrow}
\newcommand{\fpb}{0_\rightarrow}
\newcommand{\C}{\mathbb{C}}
\newcommand{\Q}{\mathbb{Q}}
\newcommand{\Z}{\mathbb{Z}}
\newcommand{\cO}{\mathcal{O}}

\newcommand{\Pp}{{\mathbf{P}^1}}
\newcommand{\Pt}{\mathbf{P}^2}
\newcommand{\cI}{\mathcal{I}}

\newcommand{\rarr}{\rightarrow}
\newcommand{\ch}{{\text{ch}}}

\DeclareMathOperator{\Hilb}{Hilb}
\DeclareMathOperator{\sgn}{sgn}

\DeclareMathOperator{\Aut}{Aut}

\newcommand{\bE}{\mathsf{E}}
\newcommand{\barQ}{\overline{Q}}
\newcommand{\barG}{\overline{G}}
\DeclareMathOperator{\rk}{rk}

\newtheorem{Theorem}{Theorem}
\newtheorem{Lemma}{Lemma}

\newtheorem{definition}[Lemma]{Definition}
\newtheorem{Proposition}[Lemma]{Proposition}

\begin{document}
\title{The local Donaldson-Thomas theory of curves}
\author{A. Okounkov and R. Pandharipande} 
\date{\begin{small}{\em Dedicated to the memory of Raoul Bott}\end{small}}
\maketitle

\begin{abstract}
The local Donaldson-Thomas theory of curves is solved
by localization and degeneration methods. The results complete
a triangle of equivalences relating Gromov-Witten theory, 
Donaldson-Thomas theory, and 
the quantum cohomology of the Hilbert scheme of
points of the plane.


\end{abstract}

\tableofcontents
\section{Introduction}

\subsection{Overview}
Let $X$ be a nonsingular projective variety of dimension 3 over $\com$.
Gromov-Witten theory  is defined by integration over the moduli space of
stable maps to $X$, and
Donaldson-Thomas theory
is defined by integration over the 
moduli space of ideal sheaves of $X$ \cite{dt,mnop1,mnop2,t}.

If $X$ is quasi-projective, the Gromov-Witten and 
Donaldson-Thomas theories may {not} be
well-defined.
However, if $X$ is the total space of a rank 2 
bundle 
over a nonsingular projective curve, 
\begin{equation*}
N\rightarrow C,
\end{equation*}
{\em local} Gromov-Witten 
and Donaldson-Thomas theories are defined via equivariant
residues \cite{bp,mnop1}.

The Gromov-Witten and  Donaldson-Thomas theories of $X$ {\em relative} to
a nonsingular 
surface $$S\subset X$$
are defined via moduli spaces of maps and sheaves with boundary 
conditions along $S$. See \cite{EGH,IP,LR,jl1,jl2,mnop2} for
various treatments of the subject.

If $X$ is the total space of a rank 2 bundle,
a natural set of 
surfaces is determined by the fibers of $N$ over points of $C$. The data
\begin{equation}\label{gret}
(N\rarr C,\ p_1, \ldots,p_r)
\end{equation}
determine relative local Gromov-Witten and Donaldson-Thomas theories 
of $N$ over $C$. 

The {\em local theory of curves} refers to {all} 
relative local Gromov-Witten and Donaldson-Thomas
theories specified by data of type \eqref{gret}.
The main result of the paper is a proof of the Gromov-Witten/Donaldson-Thomas
correspondence for the local theory of curves.

The local Gromov-Witten theory of curves is determined in \cite{bp}. Our focus
here is on the local Donaldson-Thomas theory of curves.
The paper has at least three  motivations:
\begin{enumerate}
\item[(i)] local curves provide a rich class of non-toric
examples of the GW/DT correspondence,
\item[(ii)] the proof yields a computation of the 1-legged equivariant vertex,
\item[(iii)] the correspondence for local curves 
will likely play a basic role 
 in the proof of the GW/DT correspondence for all 3-folds.
\end{enumerate}
The 1-legged equivariant vertex (ii) contains the Calabi-Yau 1-legged vertex
governed by the Gopakumar-Mari\~no-Vafa formula \cite{llz,mv,opgmv} as a special case.
The computation of the 1-legged equivariant vertex may be viewed as a
Hodge integral result on the Gromov-Witten side or a 
vertex measure result on the Donaldson-Thomas side.

The main foundational reference for  Donaldson-Thomas theory in algebraic
geometry is \cite{t}, see also \cite{mnop1,mp2}.
The required foundational development of
relative Donaldson-Thomas theory has yet
to be written. We will assume two basic properties of
relative Donaldson-Thomas theory: existence and degeneration  \cite{mnop2}.

\subsection{Definitions}
\subsubsection{Ideal sheaves}
Let $C$ be a nonsingular projective curve, and let $N$ be the
total space of a rank 2 bundle,
$$N\rightarrow C.$$
Let $I_{n} (N,d)$ denote the moduli space of of ideal sheaves
\[
0 \rightarrow I_{Z}\rightarrow  \cO _{N}\rightarrow \cO_Z \rightarrow 0
\]
of proper subschemes $Z\subset N$ of degree $d$ and Euler characteristic 
$$n=\chi (\cO _{Z}).$$
The {\em degree} of $Z$ is simply the length of the intersection
$$Z\cap N_p,$$
where $p\in C$ is a generic point.

\subsubsection{Partition functions}
If $N$ {decomposes} as a direct sum of
line bundles,
\begin{equation}\label{ffg}
N=L_1\oplus L_2,
\end{equation}
the splitting determines a scaling action of a 2-dimensional torus $T$ on $N$.
The {\em level} of the splitting is the pair of integers
$(k_1,k_2)$ where,
$$k_i= {\text {deg}}(L_i).$$
Of course, the scaling action and the level depend upon the splitting \eqref{ffg}.

Since
every $T$-fixed proper subscheme $Z$ must be supported on the
zero section of $N$,
the $T$-fixed point set, 
$$I_n(N,d)^T \subset I_n(N,d),$$ is proper.

The moduli space $I_n(N,d)$ carries a $T$-equivariant perfect obstruction 
theory obtained from (traceless) $\text{Ext}_0(I,I)$, see
\cite{t}. Though $N$ is
quasi-projective, 
$\text{Ext}_0(I,I)$ is well-defined since the associated
quotient scheme $Z\subset N$ is proper. Alternatively, for any
$T$-equivariant compactification,
$$N\subset \overline{N},$$
the obstruction theory on 
$$I_n(N,d)\subset I_n(\overline{N},d)$$
is obtained by 
restriction.

The definition of the local Donaldson-Thomas invariants of $N$ follows
the 
strategy of \cite{bp,mnop1}. We define ${\mathsf Z}(N)_{d}$ formally by:
\begin{equation}\label{eqn: moral defn of Zdt}
{\mathsf Z}(N)_{d}\text{ ``}=\text{'' }\sum _{n\in \znums
}q^{n} \int _{[I_{n} (N,d ) ]^{vir}}1 .
\end{equation}
The variable $q$  indexes the Euler characteristic $n$.
The integral on
the right of \eqref{eqn: moral defn of Zdt} is well-defined
by the virtual localization formula as an equivariant residue.

\begin{definition} The partition function for
the degree $d$ local Donaldson-Thomas invariants of $N$ is defined by:
\begin{equation}
\label{defn: definition of Zdt}
{\mathsf Z}(N)_{d}  =\sum _{n\in \znums }q^{n}
\int _{[I _{n} (N,d )^{T}]^{vir}}\frac{1}{e
(\operatorname{Norm}^{vir})}.
\end{equation}
\end{definition}

The $T$-fixed
part of the perfect obstruction theory for $I_n(N,d)$ induces a
perfect obstruction theory for $I_n (N,d )^{T}$ and hence a virtual
class \cite{gp,mnop1}. 
The 
equivariant virtual normal bundle
of the embedding, $$I_n (N,d)^{T}\subset I_n (N,d),$$ is
 $\operatorname{Norm}^{vir}$  with
equivariant Euler class
$e(\operatorname{Norm}^{vir})$. The integral in \eqref{defn: definition of Zdt}
denotes equivariant push-forward.

Following the notation of \cite{mnop1,mnop2}, 
${\mathsf Z}(N)_{d}$ is {\em unprimed} since
 the degree $0$ contributions have not yet been removed. Since a geometrical
method of removing the degree 0 contribution from 
Donaldson-Thomas theory is not available, a formal
method is followed.

\begin{definition}\label{defn: reduced DT partition function Zdt'}
The reduced partition function ${\mathsf Z}'_{DT}(N)_{d}$ for the degree $d
$ local Donaldson-Thomas invariants of $N$ is defined by:
\begin{equation*}
{\mathsf Z}'(N)_{d}  = \frac{{\mathsf Z}(N)_{d} }
{{\mathsf Z}(N)_0}. 
\end{equation*}
\end{definition}

The residues defined by the localization formula take values in the localized 
$T$-equivariant cohomology of
a point,
$$\int _{[I _{n} (N,d )^{T}]^{vir}}\frac{1}{e
(\operatorname{Norm}^{vir})} \in {\mathbb{Q}}(t_1,t_2).$$
Here, $t_1,t_2$ are the weights of the standard representations
of the factors of $T$.

If $N$ is an {\em indecomposable} rank 2 vector bundle, 
the total space of $N$ admits a scaling action of 
a $1$-dimensional algebraic torus. The local Donaldson-Thomas theory of $N$
can be defined as above with respect to the 1-dimensional scaling torus.
However, since
every indecomposable bundle $N$ is deformation equivalent
to a split bundle over $C$, the indecomposable case is
recovered from the split case via restriction to
the diagonal torus. 

In our study of the local Donaldson-Thomas
theory of $N$, we will {\em always} assume a splitting \eqref{ffg} of $N$  and
an action of a 2-dimensional scaling torus $T$.

\subsubsection{Relative geometry}
\label{relgeom}
The fiber of $N$ over a point $p\in C$ determines a $T$-equivariant
divisor 
$$N_p \subset N$$
isomorphic to $\com^2$ with the standard $T$-action.
We will consider the local theories of $N$
relative to the divisor
$$S= \bigcup_{i=1}^r N_{p_i} \subset N$$
determined by the fibers over $p_1,\ldots,p_r\in C$.

Let $I_n(N/S,d)$ denote the relative moduli space of ideal sheaves, see 
\cite{jl2, mnop2} for a discussion. The construction of $I_n(N/S,d)$, and
the existence of a
canonical $T$-equivariant $\text{Ext}_0(I,I)$ obstruction
theory will be assumed here.

For each $p_i$, let $\eta^i$
be a partition of $d$ weighted
by the equivariant Chow ring, 
$$A_T^*(N_{p_i},{\mathbb Q})\stackrel{\sim}{=} {\mathbb Q}[t_1,t_2],$$
of the fiber $N_{p_i}$.
By Nakajima's construction,
a weighted partition $\eta^i$ determines a $T$-equivariant class 
$$C_{\eta^i} \in A_T^*(\text{Hilb}(N_{p_i},d), \mathbb{Q})$$
in the
Chow ring of the Hilbert scheme of points.
In Donaldson-Thomas theory, the weighted partition $\eta^i$
specifies relative
conditions via the boundary map
$$\epsilon_i: I_n(N/S,d)\rarr \text{Hilb}(N_{p_i},d).$$

An element $\eta\in {\mathcal P}(d)$ of the set of 
partitions of $d$ may be 
viewed as a 
weighted partition with all weights set to the identity class
$1\in H^*_T(N_{p_i},{\mathbb Q})$.
The Nakajima basis of $A_T^*(\text{Hilb}(N_{p_i},d), \mathbb{Q})$ consists of 
identity weighted partitions indexed by ${\mathcal P}(d)$. 
The $T$-equivariant intersection pairing in the Nakajima basis is:
$$\int_{\text{Hilb}(N_{p_i},d)} C_\mu \cup C_\nu =
\frac{1}{(t_1t_2)^{\ell(\mu)}} 
\frac{(-1)^{d-\ell(\mu)}} 
{{\mathfrak{z}}(\mu)}\ {\delta_{\mu,\nu}},$$
where
$${\mathfrak z}(\mu) =  \prod_{i=1}^{\ell(\mu)} \mu_i \cdot 
|\text{Aut}(\mu)|.$$
The notation $\eta([0])$ will be used to set all
weights to $[0]\in A^*_T(N_{p_i},{\mathbb Q} )$.
Since
$$[0]= t_1t_2 \in A^*_T(N_{p_i}, {\mathbb Q} ),$$
the weight choice has only a mild effect.

Following the
 notation of \cite{bp,mnop2},
the relative local Donaldson-Thomas partition function,
\begin{equation*}
{\mathsf Z}(N/S)_{d,\eta^1,\dots,\eta^r}  =\sum _{n\in \znums }q^{n}
\int _{[I _{n} (N/S,d)^{T}]^{vir}}\frac{\prod_{i=1}^r \epsilon_i^*(C_{\eta^i})
}{e
(\operatorname{Norm}^{vir})},
\end{equation*}
is well-defined. 
Let
\begin{equation*}
{\mathsf Z}'(N/S)_{d,\eta^1,\dots,\eta^r}  = \frac{{\mathsf Z}(N/S)_{d,\eta^1,\dots,\eta^r } }
{{\mathsf Z}(N/S)_{0}}. 
\end{equation*}
denote the reduced relative partition function.

\subsection{Degeneration}
Simpler notation for the partition functions will often be used. If
$N$ is 
split with level $(k_1,k_2)$ over a genus $g$ base curve $C$, then let
$${\mathsf Z}(g|k_1,k_2)_{d,\eta^1,\dots,\eta^r}=
{\mathsf Z}(N/S)_{d,\eta^1,\dots,\eta^r},$$
$${\mathsf Z}'(g|k_1,k_2)_{d,\eta^1,\dots,\eta^r}=
{\mathsf Z}'(N/S)_{d,\eta^1,\dots,\eta^r}.$$
For many formulas, the $q$-shifted function,
$${\mathsf {DT}}(g|k_1,k_2)_{d,\eta^1,\dots,\eta^r}= q^{-d(1-g)} 
{\mathsf Z}'(g|k_1,k_2)_{d,\eta^1,\dots,\eta^r},$$
will be more convenient.
In the relative cases, the redundant degree subscript $d$ will
often be dropped.

Let $\bigtriangleup_d$ be the inverse of
the $T$-equivariant intersection form of the 
Nakajima basis of ${\text {Hilb}}(\com^2,d)$, 
\begin{equation}\label{mddf}
\bigtriangleup_d(\mu,\nu)=  {(-1)^{d-\ell(\mu)}}{(t_1t_2)^{\ell(\mu)}} 
{{\mathfrak{z}}(\mu)}\ {\delta_{\mu,\nu}}.
\end{equation}
The indices of the $\mathsf{DT}$ partition function
are raised by $\bigtriangleup_d$:
$${\mathsf {DT}}(g|k_1,k_2)_{\mu^1,\dots,\mu^s}^{\nu^1,\dots,\nu^t}= 
{\mathsf {DT}}(g|k_1,k_2)_{\mu^1,\dots,\mu^s, \nu^1,\dots,\nu^t}\prod_{i=1}^t
\bigtriangleup_d(\nu^i,\nu^i).$$
All the partition functions discussed here lie in the ring 
${\mathbb Q}(t_1,t_2)
((q))$ of Laurent series in $q$.

The degeneration formulas for the local
Donaldson-Thomas theory of curves are
conveniently expressed in terms of the ${\mathsf{DT}}$
partition functions:
$${\mathsf {DT}}(g|k_1,k_2)_{\mu^1,\dots,\mu^s}^{\nu^1,\dots,\nu^t} =
\sum_\gamma
{\mathsf{DT}}(g'|k'_1,k'_2)_{\mu^1,\dots,\mu^{s}}^\gamma  
{\mathsf{DT}}(g''|k''_1,k''_2)^{\nu^{1},\dots,\nu^t}_{\gamma}\, ,$$
where $g=g'+g''$, and $k_i=k_i'+k_i''$, and
$${\mathsf{DT}}(g|k_1,k_2)_{\mu^1,\dots,\mu^s} =
\sum_\gamma {\mathsf{DT}}(g-1|k_1,k_2)^\gamma_{\mu^1,\dots,\mu^s,
\gamma}\, ,$$
see \cite{mnop2} for a discussion.
The above degeneration formulas will
be assumed here.

\subsection{Localization}
Our 
localization formulas trace their origins to Bott's
remarkable paper \cite{Bott} with many stops along the
way.

The virtual localization formula of \cite{gp}, proven
in the abstract setting of perfect obstruction theories, applies to
the Donaldson-Thomas virtual class.
Applications in absolute Donaldson-Thomas
theory 
(without relative conditions) are treated
foundationally in \cite{mnop1}.
Applications in the relative setting follow from the
existence of a $T$-equivariant $\text{Ext}_0(I,I)$
obstruction theory.

\subsection{GW/DT correspondence}
Let $N$ be the a rank 2 bundle over 
a genus $g$ curve $C$ with splitting of level $(k_1,k_2)$.
The GW/DT correspondence for the local theory of curves consists of 
three results.

First, the relative local degree 0 series of $N$ is determined
in terms of the McMahon function,
$$M(q)=\prod_{n\geq 1} \frac{1}{(1-q^n)^n},$$
the generating series for 3-dimensional partitions.

\begin{Theorem}\label{tm1} The degree 0 Donaldson-Thomas partition function
is determined by:
\begin{equation*}
{\mathsf Z} (N/S)_0=M (-q)^{\int _{N}c_{3} (T_N[-S]\otimes K_N[S])}.
\end{equation*}
\end{Theorem}

Here, $T_N[-S]$ is the sheaf of tangent fields on $N$ with logarithmic zeros,
$K_N[S]$ is the logarithmic canonical bundle, and
 the integral in the exponent is defined via localization on $N$,
\begin{eqnarray*}
\int _{N}c_{3} (T_N[-S]\otimes K_N[S]) &  =  &
\int _{C}\frac{c_{3}(T_N[-S]\otimes K_N[S]) }{e ({N})} \\
&=& (2g-2+r)\frac{(t_1+t_2)^2}{t_1t_2}-(k_1+k_2), 
\end{eqnarray*}
where $r$ is the number of relative points on $C$.

Second, the reduced Donaldson-Thomas series are proven to
satisfy a basic rationality condition.

\begin{Theorem} \label{bbb}
 The reduced series ${\mathsf Z}'(N/S)_{d,\eta^1,\dots,\eta^r}$
is a rational function in the variables $t_1$, $t_2$, and $q$.
\end{Theorem}

Finally, the local Gromov-Witten theory and the local Donaldson-Thomas
theories of curves are proven to exactly match.

\begin{Theorem}\label{ccc}
After the change of variables $e^{iu}=-q$, 
\begin{multline*}
(-iu)^{d(2-2g+k_1+k_2)-\delta} 
\ {\mathsf Z}'_{GW}(N/S)_{d,\eta^1,\dots,\eta^r} 
=\\
(-q)^{-\frac{d}{2}(2-2g+k_1+k_2)}
\ {\mathsf Z}'_{DT}(N/S)_{d,\eta^1,\dots,\eta^r},
\end{multline*}
\end{Theorem}

For the Gromov-Witten side, 
we follow the definitions and notations of \cite{bp}.
In particular, 
$$\delta = \sum_{i=1}^r d-\ell(\eta^i).$$

\subsection{Method}
Theorems \ref{tm1}-\ref{ccc} are proven by solving
the local Donaldson-Thomas theory of curves. The GW/DT correspondence
is obtained by  matching
the Gromov-Witten results of \cite{bp} with the
Donaldson-Thomas results here.

The solution of the local Donaldson-Thomas theory of curves 
follows the TQFT strategy of \cite{bp}. The first step is
the determination of the level $(0,0)$ theory of $\Pp$. In the Gromov-Witten
case, integral evaluations over the moduli space of curves are required
\cite{bp}.
Parallel Donaldson-Thomas integrals are evaluated here
via connections to the quantum
cohomology of the Hilbert
scheme points of the plane. Next, the level $(-1,0)$ theory of $\Pp$
is determined by a direct calculation. Together, the results
solve the local Donaldson-Thomas theory of curves and 
prove the GW/DT correspondence.

The GW/DT correspondence for the local theory of curves has been
studied in the absolute case for the anti-diagonal action in \cite{Son}.
The correspondence for $\mathbf{P}^2$-bundles over curves
is closely related \cite{Gh,GhSon}.

\subsection{Quantum cohomology of $\text{Hilb}(\com^2,d)$} \label{heq}
For  $\lambda, \mu,\nu \in {\mathcal P}(d)$,
define the series $\langle \lambda,\mu,\nu\rangle^{\Hilb_d}$ of genus 0,
 3-pointed, $T$-equivariant Gromov-Witten invariants of 
$\text{Hilb}(\com^2,d)$
by a sum over curve degrees
$$\langle \lambda,\mu,\nu
\rangle^{\Hilb_d} = \sum_{n\geq 0} q^n
\langle \lambda,\mu,\nu
\rangle^{\Hilb_d}_{0,3,n}.$$
The insertions on the right stand for Nakajima basis elements.
See \cite{qchs} for a complete discussion of the Gromov-Witten invariants
of the Hilbert scheme $\text{Hilb}(\com^2,d)$.

The results of \cite{qchs} together with 
our calculation of the local Donaldson-Thomas theory of curves
yields a {\em Donaldson-Thomas/Hilbert
correspondence}.

\begin{Theorem}  \label{3l}
${\mathsf{DT}}(0|0,0)_{d,\lambda,\mu,\nu} 
= \langle \lambda,\mu,\nu
\rangle^{\Hilb_d}.$
\end{Theorem}

Our results complete the triangle of equivalences studied in the earlier
papers:

\begin{figure}[hbtp]\psset{unit=0.5 cm}
  \begin{center}
    \begin{pspicture}(-6,-2)(10,6)
    \psline(0,0)(2,3.464)(4,0)(0,0) 
    \rput[rt](0,0){
        \begin{minipage}[t]{4.00 cm}
          \begin{center}
            Gromov-Witten theory of $\Pp \times \C^2$ \cite{bp}
          \end{center}
        \end{minipage}}
    \rput[lt](4,0){
        \begin{minipage}[t]{3.64 cm}
          \begin{center}
             Donaldson-Thomas theory  of $\Pp \times \C^2$
          \end{center}
        \end{minipage}}
    \rput[cb](2,4.7){
        \begin{minipage}[t]{4 cm}
          \begin{center}
           Quantum cohomology of $\Hilb_d(\C^2)$ \cite{qchs}
          \end{center}
        \end{minipage}}
    \end{pspicture}
  \end{center}
\end{figure}

\noindent A fourth vertex of equivalence is obtained from the orbifold
quantum cohomology of the stack symmetric product $(\com^2)^d/\Sigma_d$,
see \cite{bg}.

\subsection{Acknowledgments} 
We thank J. Bryan, G. Farkas,
T. Graber, J. Li,  D. Maulik,  and N. Nekrasov for many  
valuable discussions. We thank J. Kock for the use of his
cobordism \LaTeX \@ macros.

Parts of the paper were written during a workshop on
algebraic geometry and topological strings at 
the Instituto Superior T\'ecnico in Lisbon in the fall of 2005.
Both authors were partially 
supported by the Packard 
foundation and the NSF.

\section{Degree 0}
\subsection{Vertex measures}
The degree 0 calculation of the local Donaldson-Thomas theory of curves
is easily obtained from the results of \cite{mnop1,mnop2}.
Let
\begin{equation}\label{fwq}
\bW(\emptyset, \emptyset,\emptyset)= 
M(-q)^{- \frac{(s_1+s_2)(s_1+s_3)(s_2+s_3)}{s_1s_2s_3}}
\end{equation}
denote the generating series of
the equivariant vertex measures of finite 3-dimensional partitions
at a 3-fold fixed point with tangent weights $s_1,s_2,s_3$.
Let 
\begin{equation}\label{fwq2}
\bW(\emptyset,\emptyset,\emptyset)_S=M(-q)^{-\frac{s_2+s_3}{s_1}}
\end{equation}
denote series of the finite vertex measures  at a 3-fold fixed point
on a relative divisor $S$ with normal weight $s_1$ and 
tangent weights $s_2,s_3$.

The evaluations \eqref{fwq} and \eqref{fwq2} are Theorem 1 and
Corollary 1 of \cite{mnop2} with an appropriate adjustment of
signs to match the 
conventions here.

\subsection{Localization}\label{zloc}
Let $S$ denote the 1-dimension 
torus acting on $\Pp$ with respective tangent weights $s$ and $-s$ at
the
fixed points $0,\infty \in \Pp$.

Let $N=L_1\oplus L_2$ be a splitting of level $(k_1,k_2)$ over $\Pp$.
The $S$-action on $\Pp$ can be lifted to $N$ with fiber weights 
$(k_1s, k_2s)$ over $0\in \Pp$ and 
fiber weights 
$(0,0)$ over $\infty \in \Pp$.
Since the scaling $T$-action on $N$ commutes with $S$,
the 3-dimensional torus,
$${\mathbf T}= S\times T,$$ 
acts on $N$.

We may calculate the degree 0 Donaldson-Thomas series 
${\mathsf Z}(0|k_1,k_2)_0$ via virtual localization with respect
to ${\mathbf T}$, see
\cite{gp,mnop1,mnop2}. 
By the evaluation of the finite vertex measure \eqref{fwq},
\begin{eqnarray*}
{\mathsf Z}(0|k_1,k_2)_0 & = 
& \Big[\bW(\emptyset,\emptyset,\emptyset)|_{{s},\ {t_1+k_1s},\ 
{t_2+k_2s}}
\cdot \bW(\emptyset,\emptyset,\emptyset)|_{-s,t_1,t_2}\Big]_{s=0} \\
& = & M(-q)^{-2\frac{(t_1+t_2)^2}{t_1t_2}- (k_1+k_2)}.
\end{eqnarray*}

The 
series ${\mathsf Z}(0|k_1,k_2)_{0,\emptyset}$ relative to $\infty \in \Pp$ is
also determined by virtual localization. Here, measure \eqref{fwq}
arises at the ${\mathbf T}$-fixed point over 0, 
and measure \eqref{fwq2} arises at the ${\mathbf T}$-fixed point
over $\infty$ of
the relative divisor:
\begin{eqnarray*}
{\mathsf Z}(0|k_1,k_2)_{0,\emptyset} 
& = & \Big[\bW(\emptyset,\emptyset,\emptyset)|_{{s},\ {t_1+k_1s},\ 
{t_2+k_2s}}
\cdot 
\bW(\emptyset,\emptyset,\emptyset)_{N_{\infty}}|_{-s,t_1,t_2}\Big]_{s=0} \\
& = & M(-q)^{-\frac{(t_1+t_2)^2}{t_1t_2}- (k_1+k_2)}.
\end{eqnarray*}

By the above evaluations, 
the proof of  Theorem \ref{tm1} is complete for the cases
${\mathsf Z}(0|k_1,k_2)_0$ and ${\mathsf Z}(0|k_1,k_2)_{0,\emptyset}$.

\subsection{Degeneration}
The degeneration formulas for the degree 0 theory take the following
two forms:
$${\mathsf Z}(g|k_1,k_2)_{0,\underbrace{\emptyset,\dots,\emptyset}_{r}} =
{\mathsf Z}(g'|k'_1,k'_2)_{0,\underbrace{\emptyset,\dots,\emptyset}_{r'},
\emptyset} \cdot 
{\mathsf Z}(g''|k''_1,k''_2)_{0,\emptyset,\underbrace{\emptyset,\dots,
\emptyset}_{r''}},$$
where $g=g'+g''$, $k_i=k_i'+k_i''$,  $r=r'+r''$, and
$${\mathsf Z}(g|k_1,k_2)_{0,\underbrace{\emptyset,\dots,\emptyset}_{r}} =
{\mathsf Z}(g-1|k_1,k_2)_{0,\underbrace{\emptyset,\dots,\emptyset}_{r},
\emptyset,\emptyset}.$$
The first degeneration formula yields a recursive equation for $r\geq 2$:
$${\mathsf Z}(0|k_1,k_2)_{0,\underbrace{\emptyset,\dots,\emptyset}_{r}} = 
\frac{{\mathsf Z}(0|k_1,k_2)_{0,\underbrace
{\emptyset,\dots,\emptyset}_{r-1}}}{{\mathsf Z}
(0|0,0)_{0,\emptyset}}.$$
From the $r=0,1$ cases, we deduce
$${\mathsf Z}(0|k_1,k_2)_{0,\underbrace{\emptyset,\dots,\emptyset}_{r}} = 
M(-q)^{(-2+r)\frac{(t_1+t_2)^2}{t_1t_2}- (k_1+k_2)}.$$
Finally, Theorem \ref{tm1} is obtained for $g>0$ by applications
of the 
second degeneration formula.
\qed

\section{TQFT}

\label{tqft}
\subsection{Overview}
The degeneration structure of local Donaldson-Thomas theory of curves
is most concisely formulated as a functor of tensor categories,
\[
\mathbf{DT} (-):2\mathbf{Cob}^{L_{1},L_{2}}\to R\mathbf{mod}.
\]
Our treatment here exactly
follows the TQFT construction in \cite{bp} for the local Gromov-Witten theory
of curves. A more detailed discussion can be found there.

\subsection{ $2\mathbf{Cob}^{L_{1},L_{2}}$}
The objects of the category $2\mathbf{Cob}^{L_{1},L_{2}}$ 
are compact oriented
1-manifolds.
A morphism in $2\mathbf{Cob}^{L_{1},L_{2}},$ 
$$Y_{1} \rightarrow Y_{2},$$ 
is an equivalence class of triples $(W,L_{1},L_{2})$ where $W$
is an oriented cobordism from $Y_{1}$ to $Y_{2}$ and $L_{1},L_{2}$ are
complex line bundles on $W$, trivialized on $\partial W$. The triples
$(W,L_{1},L_{2})$
and  $(W',L'_{1},L'_{2})$ are equivalent if there exists a boundary preserving
oriented diffeomorphism, $$f:W\to W',$$ and bundle isomorphisms $$L_{i}\cong
f^{*}L'_{i}.$$ Composition is given by concatenation of the cobordisms and
gluing of the bundles along the concatenation using the trivializations.

The isomorphism class of $L_{i}$ is determined by the Euler class $$e
(L_{i})\in H^{2} (W,\partial W),$$ which assigns an integer to each
component of $W$. For a connected cobordism $W$, we refer to the pair of
integers $(k_{1},k_{2})$, determined by the Euler classes of $L_{1}$ and
$L_{2}$, as the \emph{level}. Under concatenation, the
levels simply add. For example:

\begin{center}
\begin{texdraw}\setunitscale 1.3
\SFA 
\move (120 0)
\boundaryFA \multNB \comultNB
\move (120 30)\boundaryFA 
\move (200 0)\boundaryFA 
\move (200 30)\boundaryFA 
\htext (100 15) {$=$}
\htext (15 -5) {$\scriptstyle{(-3,1)}$}
\htext (15 30) {$\scriptstyle{(2,0)}$}
\htext (65 0) {$\scriptstyle{(7,-3)}$}
\htext (65 35) {$\scriptstyle{(-4,3)}$}
\htext (160 15) {$\scriptstyle{(2,1)}$}
\end{texdraw}
\end{center}

The empty manifold is a distinguished object in 
$2\mathbf{Cob}^{L_{1},L_{2}}$. A morphism in
$2\mathbf{Cob}^{L_{1},L_{2}}$ from the empty manifold to itself is given by
a compact, oriented, closed 2-manifold $X$ together with a pair of complex
line bundles $L_{1}\oplus L_{2}\to X$.

The category $2\mathbf{Cob}^{L_{1},L_{2}}$ is generated by the
following finite set of morphisms \cite{bp}:
\begin{center}
\begin{texdraw}\setunitscale 1.0
\unitFA
\move (30 0)
\comultFA 
\move (90 -15)
\multFA 
\move (150 0)
\counitFA 
\move (180 0)
\idFA 
\move (240 -15)
\twistFA 
\htext (-10 15) {$\scriptstyle{(0,0)}$}
\htext (45 0) {$\scriptstyle{(0,0)}$}
\htext (115 0) {$\scriptstyle{(0,0)}$}
\htext (155 15) {$\scriptstyle{(0,0)}$}
\htext (200 0) {$\scriptstyle{(0,0)}$}
\htext (245 30) {$\scriptstyle{(0,0)}$}
\htext (275 30) {$\scriptstyle{(0,0)}$}
\end{texdraw}
\end{center}

\begin{center}
\begin{texdraw}\setunitscale 1.0
\unitFA 
\move (60 0)
\unitFA 
\move (120 0)
\unitFA 
\move (180 0)
\unitFA 
\htext (-10 15) {$\scriptstyle{(0,1)}$}
\htext (50 15) {$\scriptstyle{(1,0)}$}
\htext (110 15) {$\scriptstyle{(0,-1)}$}
\htext (170 15) {$\scriptstyle{(-1,0)}$}
\end{texdraw}.
\end{center}

\subsection{The functor $\mathbf{DT}(-)$} 

Let $R$ be the ring of Laurent series 
in $q$ with coefficients given by rational
functions in $t_{1}$ and $t_{2}$,
\[
R=\mathbb{Q} (t_{1},t _{2})((q)).
\]
The collection of partition functions $\mathsf{DT} (g| k_{1},k_{2})_{\lambda
^{1},\dots, \lambda ^{r}}$ of degree $d$
gives rise to a functor
\[
\mathbf{DT} (-):2\mathbf{Cob}^{L_1,L_2}\to R\mathbf{mod}
\]
as follows. Define 
\[
\mathbf{DT} (S^{1})=H=\bigoplus _{\lambda \vdash d}Re_{\lambda }
\]
to be the free $R$-module with basis $\{e_{\lambda } \}_{\lambda \vdash d}$
labelled by partitions of $d$, and let
\[
\mathbf{DT} 
\left(S^{1} \coprod\dots \coprod S^{1} \right)=H\otimes \dots \otimes H.
\]

Let $W_{s}^{t} (g | k_{1},k_{2})$ be the connected genus $g$ cobordism from a
disjoint union of $s$ circles to a disjoint union of $t$ circles, equipped
with lines bundles $L_{1}$ and $L_{2}$ of level $(k_{1},k_{2})$. 
We define the $R$-module homomorphism
\[
\mathbf{DT} \left(W_{s}^{t} (g| k_{1},k_{2})\right):H^{\otimes s}\to H^{\otimes t}
\]
by
\[
e_{\eta ^{1}}\otimes \dots \otimes e_{\eta ^{s}}\mapsto \sum _{\mu
^{1}\dots \mu ^{t}\vdash d} \mathsf{DT} (g|k_{1},k_{2})_{\eta ^{1},\dots, \eta
^{s}}^{\mu ^{1},\dots, \mu ^{t}}e_{\mu ^{1}}\otimes \dots \otimes e_{\mu
^{t}}.
\]
We extend the definition of $\mathbf{DT} (-)$ to disconnected cobordisms
by tensor product:
\[
\mathbf{DT} \left(W[{1}]\coprod \dots \coprod W[{n}]\right)=\mathbf{DT}
\left(W[{1}])\otimes \dots \otimes \mathbf{DT} (W[{n}]\right).
\]
\begin{Proposition}\label{tqq}
$\mathbf{DT} (-): 2\mathbf{Cob}^{L_{1},L_{2}}\to
R\mathbf{mod}$ is a well-defined functor.
\end{Proposition}

\begin{proof} The degeneration formula of Donaldson-Thomas theory \cite{mnop2}
implies  the following compatibility:
\[
\mathbf{DT} \big((W,L_{1},L_{2})\circ
(W',L'_{1},L'_{2})\big)=\mathbf{DT}(W,L_{1},L_{2})\circ \mathbf{DT}
(W',L'_{1},L'_{2}).
\]
We must also prove $\mathbf{DT} (-)$ takes identity morphisms
to identity morphisms. Since the tube
$W_{1}^{1} (0|0,0)$ is the identity morphism
from $S^{1}$ to $S^{1}$ in $2\mathbf{Cob}^{L_{1},L_{2}}$, we require
\begin{equation}\label{jjww}
{\mathsf{DT}} (0| 0,0)_{\mu }^{\nu }=\delta _{\mu}^{\nu }.
\end{equation}
Equation \eqref{jjww} is proven in Lemma \ref{feww} below.
\end{proof}

\begin{Lemma}  \label{feww}
$
{\mathsf{DT}} (0| 0,0)_{\mu }^{\nu }=\delta _{\mu}^{\nu }.$
\end{Lemma}

\begin{proof}
Let $N=\cO_\Pp\oplus \cO_\Pp$ be the trivial bundle with level
$(0,0)$ splitting on $\Pp$. 
The moduli space $I_d(N/N_0\cup N_\infty,d)$ is
isomorphic to the Hilbert scheme $\text{Hilb}(\com^2,d)$.
The $q$-constant terms of ${\mathsf{DT}}(0|0,0)_{\mu,\nu}$
are therefore determined by the  
intersection form in the 
Nakajima basis:
$${\mathsf{DT}}(0|0,0)_{\mu,\nu}=
\frac{1}{   (t_1 t_2)^{\ell(\mu)}      }
  \frac{(-1)^{d-\ell(\mu)}}
{\mathfrak{z}(\mu)} \ \delta_{\mu,\nu}.$$
Hence, the Lemma is proven for $q$-constant terms.

The degeneration formula in Donaldson-Thomas theory
yields the following factorization
\begin{equation*}
  \sum_\nu   {\mathsf{DT}} (0| 0,0)_{\mu }^{\nu }
{\mathsf{DT}} (0| 0,0)_{\nu }^{\rho } = {\mathsf{DT}} (0| 0,0)_{\mu }^{\rho }.
\end{equation*}
Since the matrix
${\mathsf{DT}} (0| 0,0)_{\mu }^{\nu }$ is invertible by the
$q$-constant analysis,
${\mathsf{DT}} (0| 0,0)_{\mu }^{\nu }$ must be the identity
matrix (with no $q$-dependence).
\end{proof}

\section{Vanishing}
\subsection{Summary}
Consider the local Donaldson-Thomas theory of  
level $(0,0)$ on a
nonsingular curve $C$.
Let $$T=\com^* \times \com^*$$ be the 2-dimensional scaling torus.
Let $\TT$ denote the 1-dimensional  anti-diagonal subtorus,
$$\TT = \{ (\xi, \xi^{-1})\   | \ \xi \in \com^*\} \subset T.$$
The $\TT$-equivariant Donaldson-Thomas invariants are obtained
from the $T$-equivariant invariants by the substitution
$$t_1=t, \ \ t_2=-t,$$
where $t$ is the weight of the standard representation of $\TT$.
We prove vanishing results for the $\TT$-equivariant
Donaldson-Thomas invariants of level $(0,0)$.

\subsection{Descendent insertions}
Let $N$ be a split rank 2 bundle on a nonsingular curve $C$ with
a scaling $T$-action.
We will consider the local 
Donaldson-Thomas theory of $N$ with descendent insertions. We review
the definitions of \cite{mnop2}.

The moduli space $I_n(N,d)$ is canonically isomorphic to the
Hilbert scheme of curves of $N$, see \cite{mp2}. 
Let $\pi_1$ and $\pi_2$ denote the projections to the respective
factors of $I_n(N,d) \times N$.
Consider the universal ideal sheaf $\mathfrak{I}$,
$$\mathfrak{I} \rarr I_n(N,d) \times N.$$
Since $\mathfrak{I}$ is $\pi_1$-flat 
and $N$ is nonsingular, a finite resolution
of $\mathfrak{I}$ by locally free sheaves on $I_n(N,d) \times N$
exists. Hence, the Chern classes of $\mathfrak{I}$ are well-defined.

For $\gamma\in A^l_T(N,\mathbb{Q})$, 
let $\ch_{k+2}(\gamma)$ denote the following operation on the 
Chow homology of $I_n(N,d)$:
$$\ch_{k+2}(\gamma):A^T_*(I_n(N,d),\mathbb{Q}) 
\rarr A^T_{*-k+1-l}(I_n(N,d),\mathbb{Q}),$$
$$\ch_{k+2}(\gamma)\big(\xi\big)=
 \pi_{1*}\big(   \ch_{k+2}(\mathfrak{I})
\cdot \pi_2^*(\gamma) \cap      \pi_1^*(\xi)\big).$$
Though $\pi_1$ is {\em not} proper, 
the $T$-equivariant push-forward $\pi_{1*}$ is well-defined by localization.

Descendent fields in Donaldson-Thomas theory, denoted by
${\sigma}_k(\gamma)$,  correspond to the operations $\ch_{k+2}(\gamma)$. 
The $T$-equivariant descendent invariants are defined by
\begin{equation}\label{dxcc2}
\Big\langle {\sigma}_{k_1}(\gamma_{l_1}) \cdots
{\sigma}_{k_r}(\gamma_{l_r})\Big\rangle^{N}_{n,d} 
= \int_{[{I}_n(N,d)^T]^{vir}} 
\frac{\prod_{i=1}^r \ch_{k_i+2}(\gamma_{l_i})}{e(\text{Norm}^{vir})},
\end{equation}
where the latter integral is the push-forward to a point of
the class
$$\ch_{k_1+2}(\gamma_{l_1}) \ \circ\ \cdots\ \circ
\ch_{k_r+2}(\gamma_{l_r})\Big( \frac{[I_n(N,d)^T]^{vir}}
{e(\text{Norm}^{vir})}\Big).$$
The descendent invariants of $N$ may be viewed equivalently
as equivariant residues:
$$
\Big\langle {\sigma}_{k_1}(\gamma_{l_1}) \cdots
{\sigma}_{k_r}(\gamma_{l_r})\Big\rangle^{N}_{n,d} 
= \text{Res}_{I_n(N,d)^T} \left[ \int_{[{I}_n(N,d)]^{vir}} 
{\prod_{i=1}^r \ch_{k_i+2}(\gamma_{l_i})}\right].
$$

The definition of $T$-equivariant descendent invariants in relative
Donaldson-Thomas theory of $N$ is identical.
The boundary condition over a relative point $p_i\in C$ is
determined by a partition $\eta^i$ weighted by $H^*_T(N_{p_i}, {\mathbb Q})$.

Brackets with relative conditions on the right will often be used.
For example,
\begin{equation}\label{cqww}
\Big\langle {\sigma}_{k_1}(\gamma_{l_1}) \cdots
{\sigma}_{k_r}(\gamma_{l_r})\ \Big| \ \nu^1, \ldots, \nu^s
\Big\rangle^{N}_{n,d} 
\end{equation}
denotes a descendent invariant relative to $s$ points of $C$.

\subsection{Brackets}\label{bcon}
Efficient bracket notation for Donaldson-Thomas invariants will be
used throughout the paper.

For absolute brackets (without relative conditions), the degree subscript $d$
is always required. If the Euler characteristic subscript $n$ is
omitted, a sum is signified,
$$\Big\langle \prod_i {\sigma}_{k_i}(\gamma_{l_i})
\Big\rangle^{N}_{d} = \sum_{n} q^n 
\Big\langle \prod_i {\sigma}_{k_i}(\gamma_{l_i})
\Big\rangle^{N}_{n,d}.$$

If a relative condition occurs in a bracket, the degree subscript is
redundant and therefore may be omitted,
$$
\Big\langle \prod_i {\sigma}_{k_i}(\gamma_{l_i})
\ \Big| \ \nu
\Big\rangle^{N}_{n} =
 \Big\langle \prod_i {\sigma}_{k_i}(\gamma_{l_i})
\ \Big| \ \nu
\Big\rangle^{N}_{n, |\nu|}.
$$
If all subscripts in a relative bracket are omitted, a sum is
signified,
$$
\Big\langle \prod_i {\sigma}_{k_i}(\gamma_{l_i})
\ \Big| \ \nu
\Big\rangle^{N} = \sum_{n} q^n
 \Big\langle \prod_i {\sigma}_{k_i}(\gamma_{l_i})
\ \Big| \ \nu
\Big\rangle^{N}_{n, |\nu|},
$$
as in the absolute case.

Most of our Donaldson-Thomas computations will concern
the local theory of $\Pp$. If the superscript $N$ is
replaced by a level $(k_1,k_2)$, 
the theory of $\Pp$ is signified,
$$
\Big\langle \prod_i {\sigma}_{k_i}(\gamma_{l_i})
\ \Big| \ \nu
\Big\rangle^{(m_1,m_2)} = 
\Big\langle \prod_i {\sigma}_{k_i}(\gamma_{l_i})
\ \Big| \ \nu
\Big\rangle^{\cO_\Pp(m_1)\oplus \cO_\Pp(m_2)}.$$
If the superscript is omitted altogether, then
the level $(0,0)$ theory of $\Pp$ is signified,
$$
\Big\langle \prod_i {\sigma}_{k_i}(\gamma_{l_i})
\ \Big| \ \nu
\Big\rangle = 
\Big\langle \prod_i {\sigma}_{k_i}(\gamma_{l_i})
\ \Big| \ \nu
\Big\rangle^{(0,0)}.$$
Of course, redundant labels may be kept in various
formulas for emphasis.

\subsection{Restriction to $\TT$}

Consider the level $(0,0)$ theory on a nonsingular curve $C$.
Since the $\TT$-fixed locus of $I_n(\cO_C\oplus\cO_C,d)$ is proper,
the $\TT$-equivariant descendent invariants are well-defined
by residues. 
The restriction of the ${T}$-equivariant descendent 
invariants 
to $\TT$  yield the $\TT$-equivariant descendent invariants.
The restriction to $\TT$ can also be seen to be well-defined
by the following more precise result.

\begin{Lemma} The relative
$T$-equivariant descendent invariants
of level $(0,0)$ on $C$ take values in the subring
$${\mathbb Q}[t_1,t_2]_{(t_1t_2)}\subset {\mathbb Q}(t_1,t_2).$$
\end{Lemma}

\begin{proof} 
Let $N=\cO_C \oplus \cO_C$.
As before, let
$${\mathfrak J} \rarr {\mathcal N}$$
denote the universal ideal sheaf over the universal total space
$${\mathcal N} \rarr I_n(N/S,d).$$ 
Since $N=C \times \com^2$, there is proper morphism
$$\rho:{\mathcal N} \rarr I_n(N/S,d) \times \com^2.$$
Moreover,
$\rho_*(\cO_{\mathcal N}/{\mathfrak J})$ is flat family over
$I_n(N/S,d)$ of
torsion sheaves of $\com^2$ of length $n+d g_C$.
There is an associated
morphism of Hilbert-Chow type,
$$\iota:I_n(N/S,d) \rarr \text{Sym}^{n+dg_C}(\com^2).$$
A $T$-equivariant, proper
morphism,
$$\iota': \text{Sym}^{n+d g_C}(\com^2) \rarr \oplus_{1}^{n+dg_C}
\com^2,$$
is obtained via the higher moments,
$$\iota'\Big( \ \{(x_i,y_i)\} \ \Big) = \Big(\sum_i x_i, \sum_i y_i\Big) \oplus 
\Big(\sum_i x^2_i, \sum_i y^2_i\Big) \oplus \cdots \oplus 
\Big(\sum_i x^n_i, \sum_i y^n_i\Big).$$
Let $j=\iota'\circ \iota$.

Since $j$ is a $T$-equivariant, proper morphism,
there is $T$-equivariant push-forward
$$j_*: A^T_*(I_n(N/S,d), {\mathbb Q}) \rarr
A^T_*(   \text{Sym}^{n+d\cdot g(C)}(\com^2)    , {\mathbb Q}).$$
Descendent invariants
are defined via the $T$-equivariant residue of
$$\left(\prod_i \text{ch}_{k_i}(\gamma_{l_i}) \cup \epsilon_{\text{rel}}\right) \ \cap
[I_n(N/S,d)]^{vir} \ \in A^T_*(I_n(N/S,d), {\mathbb Q}),$$
where $\epsilon_{\text{rel}}$ denotes the relative conditions.
We may instead calculate the $T$-equivariant residue of
$$j_*\left( \left( 
\prod_i \text{ch}_{k_i}(\gamma_{l_i}) \cup \epsilon_{\text{rel}}\right) \ \cap
[I_n(N/S,d)]^{vir}\right) \ \in A^T_*( 
\oplus_1^{n+dg_C}(\com^2)    , {\mathbb Q}).$$
Since the space $\oplus_1^{n+d\cdot g(C)} 
\com^2$ has a unique $T$-fixed point
with tangent weights,
$$t_1,t_2,2t_1,2t_2, \ldots, (n+dg_C)t_1, (n+dg_C)t_2,$$
we conclude the descendent invariant
has only monomial poles in the variables $t_1$ and $t_2$.
\end{proof}

We denote the restriction of the $T$-equivariant
descendent invariants to the anti-diagonal subtorus 
by  an additional superscript $\pm$. For
example, the restriction of \eqref{cqww} is denoted by
$$
\Big\langle {\sigma}_{k_1}(\gamma_{l_1}) \cdots
{\sigma}_{k_r}(\gamma_{l_r}) \ \Big| \ \nu^1, \ldots, \nu^s
\Big\rangle^{N \pm}_{n,d} .$$

\subsection{Absolute $\Pp$}
The bundle $N=\cO_\Pp \oplus \cO_\Pp$ admits a natural action of
the 3-dimensional torus,
$${\mathbf T}= S \times  T,$$
via the canonical
lifting of $S$, see Section \ref{zloc}.
The $T$-equivariant descendent invariants of $N$ can
be calculated by localization on $I_n(N,d)$ with respect to
the ${\mathbf T}$-action. 

The localization of the virtual 
class  $[I_n(N,d)]^{vir}$ to the 
${\mathbf T}$-fixed points of $I_n(N,d)$
 is determined by the 
formulas of \cite{mnop1, mnop2} in terms of vertex 
and edge terms. 

Consider first the vertex terms.
Let
$\Pi(\lambda,\emptyset,\emptyset)$ be the set
of  3-dimensional partitions with  outgoing
2-dimensional partitions $\lambda$, $\emptyset$,
and  $\emptyset$.
The partitions $\pi\in \Pi(\lambda,\emptyset,\emptyset)$ are
 finite in two of the three outgoing
directions. The generating series $\bW(\lambda,\emptyset,\emptyset)$
is defined by
$$\bW(\lambda,\emptyset,\emptyset) =
\sum_{\pi\in \Pi(\lambda,\emptyset,\emptyset)} \bw(\pi) q^{|\pi|}.$$
Here, $\bw(\pi)$ is the equivariant vertex measure, and
 $|\pi|$ is the number of boxes of $\pi$ which remain
after removing the infinite outgoing cylinder \cite{mnop1,mnop2}.

\begin{Lemma}\label{xtyy}
For $\pi\in \Pi(\lambda,\emptyset,\emptyset)$ satisfying
$|\pi|>0$, the measure $\bw(\pi)|_{s,t_1,t_2}$ is divisible by $t_1+t_2$.
\end{Lemma}

\begin{proof}
The proof exactly follows the derivation of Lemma 4 of \cite{mnop2}.
We determine here
the precise positive power of $t_1+t_2$ 
dividing the measure $\bw(\pi)|_{s,t_1,t_2}$.

Let $\lambda$ and $\mu$ be two partition diagrams
satisfying 
$\lambda\supset\mu$. The difference $\lambda/\mu$
 between $\lambda$ and $\mu$ is a
\emph{skew diagram}. 
The content $c(\square)$ of a square of a partition diagram
with coordinates $(i,j)$ is defined by
$$
c(\square) = j-i \,.
$$
A \emph{rim hook} is a connected skew diagram
which does not 
contain two squares of equal content.

For any skew diagram $\lambda/\mu$ there is 
a minimal integer $r$ for which
$$
\mu = \nu_0 \subset \nu_1 \subset \dots \subset \nu_r = \lambda
$$
and each $\nu_{k+1}/\nu_k$ is a rim hook. The minimal $r$ is
the 
\emph{rank} of $\lambda/\mu$. 
The 
rank can be determined by repeatedly peeling off maximal rim hooks 
from $\lambda$. The process can be seen in a rank 4 example:

  \begin{center}
    \scalebox{0.7}{\includegraphics{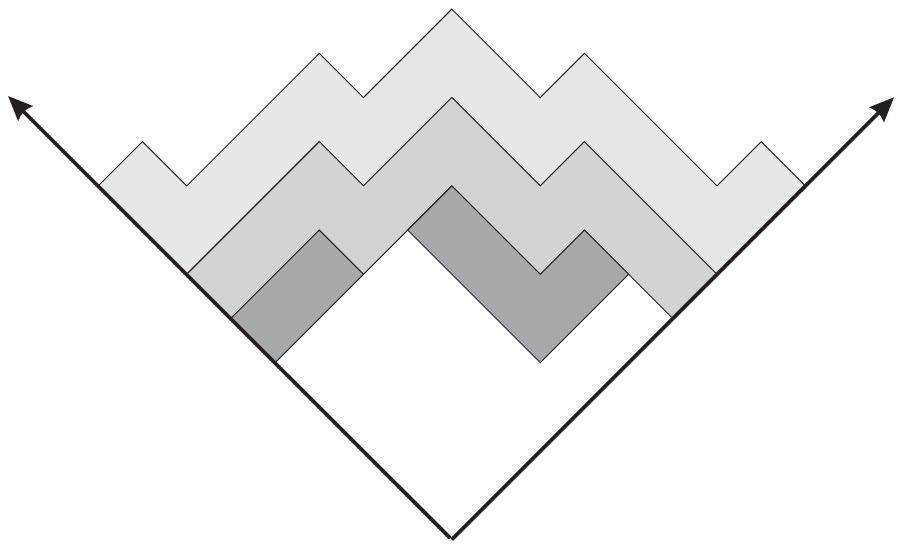}}
  \end{center}

Let $a_k$ be the number of squares in $\lambda/\mu$ of content 
$k$. The rank is determined in terms of $a_k$ by
\begin{equation}\label{sum_cont}
\rk \lambda/\mu =  \frac12 \sum_k |a_k-a_{k+1}|
\end{equation}
since
the right side of \eqref{sum_cont} receives contributions
from the beginning and end of
every rim hook. Since
each summand in \eqref{sum_cont} is either $0$ or $1$, 
each term can be squared
\begin{equation}\label{sum_cont2}
\rk \lambda/\mu =  \frac12 \sum_k (a_k-a_{k+1})^2 \,.
\end{equation}

A 3-dimensional partition with one leg of profile $\lambda$ 
can be viewed as a sequence of slices by planes 
perpendicular to the direction of the leg: 
$$
\lambda^{0} \supset \lambda^{(1)} \supset \lambda^{(2)} 
\supset \lambda^{(M)} = \dots = \lambda^{(\infty)} = \lambda \,, 
\quad M\gg 0\,. 
$$
The order of divisibility of $\bw(\pi)|_{s,t_1,t_2}$
by $t_1+t_2$ is
\begin{equation}\label{ord_van}
  \textup{ord}_{(t_1+t_2)} \, \bw(\pi)|_{s,t_1,t_2} = 
\sum_{k=0}^\infty \rk \lambda^{(k)}/\lambda^{(k+1)} \,,
\end{equation}
where, in fact, only finitely many terms are nonzero. 
Formula \eqref{ord_van} is a direct consequence 
of \eqref{sum_cont2} and the proof of Lemma 4 of  \cite{mnop2}. 
\end{proof}

The ${\mathbf T}$-fixed points of 
$I_n(N,d)$ are isolated and  correspond bijectively
to triples
$(\pi,\lambda, \pi')$
where 
$\pi,\pi' \in \Pi(\lambda,\emptyset,\emptyset)$ and
$$|\pi|+|\pi'|=n-|\lambda|.$$
The vertex partitions $\pi, \pi'$  determine the nonreduced structure of
the ${\mathbf T}$-fixed ideal over $0,\infty \in \Pp$.
The edge partition $\lambda$ determines the nonreduced structure
over $\Pp \setminus \{0,\infty\}$.

The ${\mathbf T}$-equivariant
localization of the virtual class $[I_n(N,d)]^{vir}$ to the fixed point
$(\pi,\lambda, \pi')$ is
\begin{equation}
\label{jjw}
\bw(\pi)|_{s,t_1,t_2} \cdot {E(\lambda)} 
\cdot \bw(\pi')|_{-s,t_1,t_2},
\end{equation}
where the edge terms $E(\lambda)$ 
are the (inverse) tangent $T$-weights to $\text{Hilb}(\com^2, |\lambda|)$
at the $T$-fixed point 
indexed by $\lambda$. The edge terms are 
easily seen to be prime to $t_1+t_2$.

By Lemma \ref{xtyy}, the localization \eqref{jjw} vanishes when restricted
to the anti-diagonal torus $\TT$ if either $|\pi|$ or $|\pi'|$ are
positive.

\begin{Lemma} \label{ooo}
For $n>d$, 
$$\Big\langle {\sigma}_{k_1}(\gamma_{l_1}) \cdots
{\sigma}_{k_r}(\gamma_{l_r})
\Big\rangle^{\cO_\Pp \oplus \cO_\Pp  \pm}_{n,d} = 0.$$
\end{Lemma}

\begin{proof}
The result is a consequence of the vanishing \eqref{jjw} of 
the virtual class
when localized to the ${\mathbf T}$-fixed points (and restricted to
$\TT$).
The descendent integrand plays no
role.
\end{proof}

Lemma \ref{ooo} is the first 
vanishing result for the $\TT$-equivariant Donaldson-Thomas
invariants of level $(0,0)$.

\subsection{The matrix $M_d$} \label{thmx}
Let $N=\cO_\Pp \oplus \cO_\Pp$ be the trivial bundle with
level $(0,0)$ splitting over $\Pp$. 
Let $N_0$ denote the fiber of $N$ over
$0\in \Pp$. Let
$$[N_0]\in A_T^1(N,{\mathbb Q})$$
be the associated class.

We define a matrix $M_d$ of 
descendent invariants of $N$ relative to $N_\infty$
indexed by the set ${\mathcal P}(d)$ of partitions of $d$.
For partitions $\mu,\nu\in {\mathcal P}(d)$, let
$$M_d(\mu,\nu) =  q^{-d} 
\left. \left \langle \prod_i
\sigma_{\mu_i-1}([N_0])\ \right| \ \nu([0])\right \rangle^N,$$
following the bracket conventions of Section \ref{bcon}.
The partition $\mu$ specifies a descendent insertion, and the partition $\nu$
specifies a relative condition along $N_\infty$. Since,
the minimal Euler characteristic of a degree $d$ subscheme of $N$ is $d$,
the elements of $M_d$ lie in the ring 
${\mathbb Q}[t_1,t_2]_{(t_1t_2)}[[q]]$.

Define 
the {\em length} partial order on ${\mathcal P}(d)$ by the following rule:
$\mu\geq  \mu'$ if $\ell(\mu)> \ell(\mu')$ or if $\mu=\mu'$.

\begin{Lemma} $M_d$ is 
upper triangular with respect to the length partial ordering.
\end{Lemma}

\begin{proof} 

Let $N\subset \overline{N}$ denote the $T$-equivariant compactification
over $\Pp$
defined by
$$\overline{N}={\mathbf P}(N \oplus \cO_\Pp).$$
 Here,
$T$ acts trivially on the additional $\cO_\Pp$.
We will only consider the curve classes on $\overline{N}$
obtained from $N$.

Let $\overline{N}_\infty$ denote the fiber of the compactification
over $\infty\in \Pp$.
Let 
$$\fpa, \fpb \in \overline{N}_\infty$$ be the two {\em new}
$T$-fixed points with normal
$T$-weights $$t_1-t_2, -t_2, \ \text{and} \ t_2-t_1,-t_1$$
respectively.

By the residue definition,
the $T$-equivariant Donaldson-Thomas descendent invariant
\begin{equation}\label{xde}
\left. \left \langle \prod_i \sigma_{\mu_i-1}([N_0]) 
\ \right| \ \nu([0]) \right\rangle^N_{n}
\end{equation}
occurs as a 
summand in the localization computation of the $T$-equivariant descendent
\begin{equation}\label{xcv}
\left. \left \langle \prod_i \sigma_{\mu_i-1}([\overline{N}_0]) 
 \ 
\right| \ \nu([0]) \right \rangle^{\overline{N}}_{n}
\end{equation}
for {\em any} partition $\mu$ --- not necessarily a partition of $d$.

The virtual dimension of the moduli space $I_n(\overline{N}/\overline{N}_\infty
,d)$ is $2d$. 
The integrand and relative constraints of \eqref{xcv} impose 
$$|\mu|-\ell(\mu) + d + \ell(\nu)$$
conditions.
Therefore, since $\overline{N}$ is proper, the integral \eqref{xcv} vanishes
if $$|\mu|-\ell(\mu)+ \ell(\nu)<d.$$

The ${T}$-localization formula expresses \eqref{xcv} as the
following sum of triple products of
relative local Donaldson-Thomas invariants of $N$:
\begin{multline}
\sum_{n_0+n_\uparrow+n_{\rightarrow}=n} \ 
\sum_{A_0 \cup A_\uparrow \cup {A}_{\rightarrow}
=\{1,\ldots, \ell(\mu)\}} \label{zxz}
\left \langle \left.  \prod_{i\in A_0} \sigma_{\mu_i-1}([{N}_0]) 
\ \right| \ \nu([0]) \right \rangle^{{N}}_{n_0} \\
\hspace{+125pt}
\cdot 
 \left\langle \left. \prod_{i\in A_\uparrow}\sigma_{\mu_i-1}([{N}_0]) 
 \ \right|\ \emptyset  \right\rangle^{{N}}_{n_\uparrow} \Big|_{t_1-t_2,-t_2}\\
\hspace{-90pt}
\cdot 
\left \langle \left. \prod_{i\in A_\rightarrow}
\sigma_{\mu_i-1}([{N}_0]) \ \right|
\ \emptyset \right \rangle^{{N}}_{n_\rightarrow}\Big|_{t_2-t_1,-t_1}.
\end{multline} 
The relative localization formula is applied  with
a factorization rule ---
the $T$-fixed loci do {\em not} geometrically factor.\footnote{
In the absolute case, the $T$-fixed loci of
$I_n(\overline{N},d)$
factor as a triple product.
The $T$-fixed loci relative
to a fiber,
$I_n(\overline{N}/\overline{N}_\infty,d)$
do {\em not} factor. However, a
factorization rule holds.
Factorization can be deduced from the relative localization
formula applied to disjoint unions following the
discussion of connected/ disconnected
issues in relative Gromov-Witten theory in Section 1.8 of \cite{mp}.}

By induction on $n$ and $\ell(\mu)$, we conclude
\begin{equation}\label{van7}
\left \langle\left. \prod_i \sigma_{\mu_i-1}([{N}_0]) \ 
\right| \ \nu([0]) \right \rangle^{{N}}_{n}=0
\end{equation}
if $|\mu|-\ell(\mu)+\ell(\nu)<d$.
In particular,
$$M_d(\mu,\nu)=0$$
if $\ell(\mu)> \ell(\nu)$.

If $|\mu|-\ell(\mu)+\ell(\nu)=d$, the constraints of
\eqref{xcv} impose exactly $2d$ conditions. The nonequivariant
integral
\begin{equation}\label{xcvv}
\left \langle\left. \prod_i\sigma_{\mu_i-1}([\overline{N}_0])  \ \right| 
 \ \nu_1([\gamma]), \nu_2([0]),
\ldots, \nu_{\ell(\nu)}([0])\right \rangle^{\overline{N}}_{n}
\end{equation}
is independent of
$\gamma \in \overline{N}_\infty$ --- the parts of the
relative condition over $\infty$ are written explicitly here.
After specializing $\gamma$ to a 
$T$-fixed point of $\overline{N}_\infty$, 
the invariant \eqref{xcvv} can be computed by localization with respect to $T$.

If $\gamma$ is specialized to $0\in 
\overline{N}_\infty$,  $T$-equivariant
localization expresses \eqref{xcvv}, as before, as a 
sum of triple product \eqref{zxz}.
The vanishing \eqref{van7} removes most terms.

If $\gamma$ 
is specialized to $\fpa$,
the $T$-equivariant localization formula for \eqref{xcvv} takes
a different form:
\begin{multline}
\sum_{n_0+n_\uparrow+n_{\rightarrow}=n} \
\sum_{A_0 \cup A_\uparrow \cup A_{\rightarrow}=
\{1,\ldots, \ell(\mu)\}}  \label{zxzz}
\left \langle \left. \prod_{i\in A_0} \sigma_{\mu_i-1}([{N}_0]) 
\ \right| \ \nu_2([0]),\dots,\nu_{\ell(\nu)}([0]) 
\right \rangle^{{N}}_{n_0} \\
\hspace{+150pt}
\cdot 
\left \langle \left. \prod_{i\in A_\uparrow}
\sigma_{\mu_i-1}([{N}_0]) 
\ \right| \ \nu_1([0]) \right \rangle^{{N}}_{n_\uparrow}
\Big|_{t_1-t_2,-t_2}\\
\cdot 
\left \langle \left. \prod_{i\in A_{\rightarrow}}\sigma_{\mu_i-1}([{N}_0]) 
\ \right| \ \emptyset \right \rangle^{{N}}_{n_{\rightarrow}}
\Big|_{t_2-t_1,-t_1}.
\end{multline}

By repeated use of the 
comparison of the two evaluations of \eqref{xcvv}, we find
$$\left \langle \left. \prod_{i} \sigma_{\mu_i-1}([{N}_0]) 
\ \right| \ \nu([0]) \right \rangle^{{N}}_{n}=0$$
unless there are disjoint subpartitions $\mu[i]\subset \mu$
such that
\begin{equation}\label{pppz}
|\mu[i]|-\ell(\mu[i])+1=\nu_i.
\end{equation}
If $\mu\in{\mathcal P}(d)$ and $\ell(\mu)=\ell(\mu)$, condition
\eqref{pppz} implies $\mu=\nu$.
Hence,
$$M_d(\mu,\nu)=0$$
if $\ell(\mu)=\ell(\nu)$ unless $\mu=\nu$.
\end{proof}

\begin{Lemma}  $M_d$ is invertible in the ring of matrices with  
${\mathbb Q}[t_1,t_2]_{(t_1t_2)} [[q]]$ coefficients.
\end{Lemma}
\begin{proof}
The minimal Euler characteristic of a degree $d$ subscheme of $N$ is $d$. Since
$$I_d(N,d) \stackrel {\sim}{=} \text{Hilb}(N_\infty, d),$$
the matrix of $q$-constant terms of $M_d$ is 
determined by the classical (equivariant) intersection
theory of the Hilbert scheme of points of the plane and is well-known to
be invertible.
\end{proof}

Let $M_d^{\pm}$ denote the restriction of $M_d$ to the anti-diagonal
torus. The following vanishing result holds.

\begin{Lemma}\label{ttt}
$M_d^{\pm}$ has no $q$ dependence.
\end{Lemma}

\begin{proof} Let $C_{d}$ be a matrix indexed by partitions
${\mathcal P}(d)$ with the coefficients
$$C_d(\mu,\nu)= q^{-d}
\left\langle \prod_{i}
{\sigma}_{\mu_i-1}([N_0]) \  \cdot\   \prod_{j}{\sigma}_{\nu_j-1}([N_\infty]) 
   \right\rangle^{N}_{d}.$$
The degeneration formula in Donaldson-Thomas theory,
yields the factorization
\begin{equation} \label{dsvef}
M_d \bigtriangleup_d M_d^t = C_d,
\end{equation}
where $\bigtriangleup_d$ is defined by \eqref{mddf}.
The matrix  $\bigtriangleup_d$ has
 {\em no} $q$ dependence.

Equation \eqref{dsvef} is obtained by a degeneration of the 
base 
to a reducible nodal curve:

  \begin{center}
    \scalebox{0.7}{\includegraphics{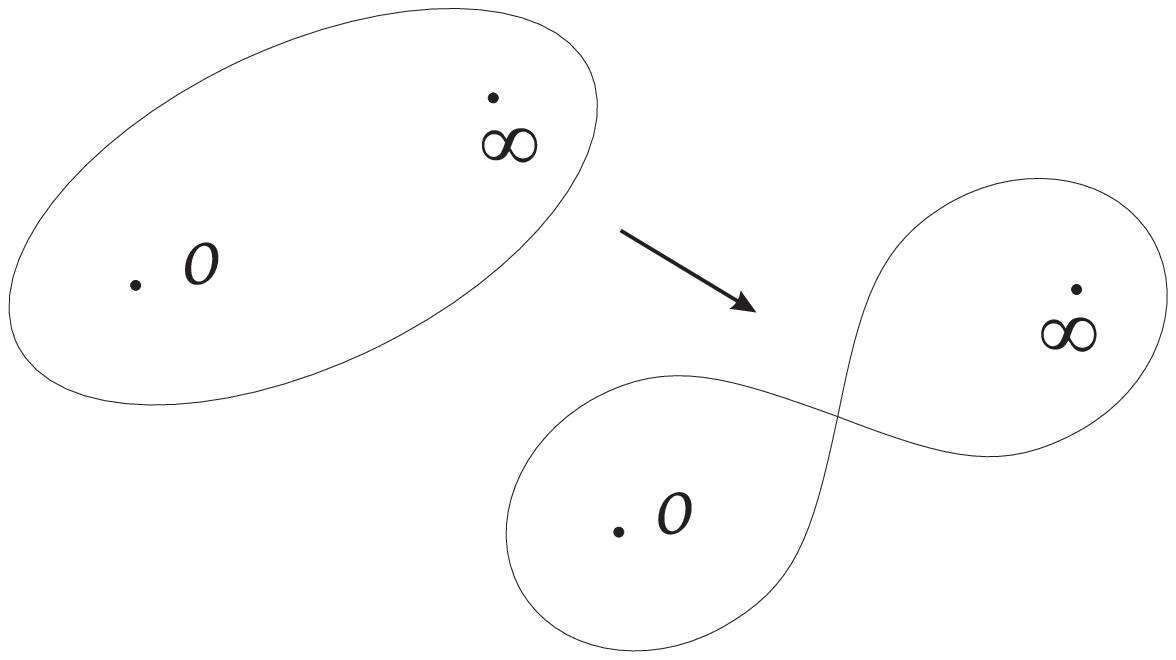}}
  \end{center}

\noindent The bundle $N$ specializes to a sum of trivial bundles on the
reducible curve. The degeneration is equivariant for the scaling torus
$T$.

The restriction $C_d^{\pm}$ has no $q$ dependence by Lemma \ref{ooo}.
The $\TT$-restriction of
\eqref{dsvef} is a Gauss decomposition of $C_d^\pm$. Uniqueness
of the Gauss decomposition implies the Lemma.
\end{proof}

Lemma \ref{ttt} can be restated as an explicit vanishing of
Donaldson-Thomas invariants. For $n>d$, 
\begin{equation} \label{vqa}
 \left 
\langle \left. \prod_i \sigma_{\mu_i-1}([N_0]) \ \right| \ 
\nu([0])\right \rangle^{\cO_\Pp\oplus \cO_\Pp\pm}_{n,d} =0,
\end{equation}
where $\mu,\nu \in \mathcal{P}(d)$.

\subsection{Degeneration}
Let $N=\cO_C \oplus \cO_C$ be the trivial bundle with level $(0,0)$ splitting
over a nonsingular, genus $g$ curve $C$. The class 
$$[N_z]\in A^1_T(N,{\mathbb Q})$$
is independent of $z\in C$. We will denote the fiber
class $[N_z]$ by $F$.

\begin{Lemma} For $n>d(1-g)$, \label{hx}
\begin{equation}\label{h23}
\Big 
\langle \prod_{i} \sigma_{k_i}(F)\ \Big| \ 
\nu^1,\ldots, \nu^s \Big \rangle^{\cO_C\oplus \cO_C\pm}_{n,d} =0.
\end{equation}
\end{Lemma}

\begin{proof}
Since the matrix $M_d^{\pm}$ is invertible and has no $q$ dependence,
the relative conditions in the Donaldson-Thomas integral \eqref{h23} can
be systematically traded for descendent insertions. The Lemma is then
equivalent to the vanishing for $n>d(1-g)$ of {\em all} absolute invariants
\begin{equation}\label{h234}
\Big 
\langle \prod_{i'} \sigma_{k'_{i'}}(F) \Big \rangle^{\cO_C\oplus \cO_C 
\pm}_{n,d} =0.
\end{equation}
After degenerating $C$ to a nodal rational curve (and again trading
relative conditions for descendent insertions), we need only prove the vanishing
\eqref{h234} in case $C$ is $\Pp$. The latter vanishing is a consequence of
Lemma \ref{ooo}.
\end{proof}

\subsection{Cotangent lines, rubber, and topological recursion}\label{clrt}
Let $N=\cO_C \oplus \cO_C$ be the trivial bundle with 
level $(0,0)$ splitting over
 a nonsingular curve $C$.
Consider the moduli space of ideal sheaves $I_n(N/N_p,d)$
relative to the fiber over $p\in C$.
The cotangent line bundle,
$${\mathbb L}_p \rarr I_n(N/N_p,d),$$
is defined by the cotangent space at the relative point
$p$ of the target curve.
The torus $T$ acts trivially on ${\mathbb L}_p$.
Let
$$\psi_p \in A^1_T(I_n(N/N_p,d), {\mathbb Q})$$
denote the first Chern class of ${\mathbb L}_p$.

The Donaldson-Thomas theory of {\em rubber} naturally arises at the
boundary of $I_n(N/N_p,d)$.
Let $R$ be a rank 2 bundle of level $(0,0)$ over $\Pp$. Let 
 $$R_0, R_\infty\subset R$$ 
denote the fibers over $0, \infty\in \Pp$.
The 1-dimensional torus $S$ acts on $R$ via the symmetries of
$\Pp$. 
Let $I_n(R/R_0\cup R_\infty,d)$ be the relative moduli space
of ideal sheaves, and let
$$I_n(R/R_0 \cup R_\infty,d)^\circ \subset I_n(R/R_0\cup R_\infty,d)$$
denote the open set with finite stabilizers for the $S$-action
and {\em no} destabilization over $\infty\in \Pp$.
The rubber moduli space,
$${I_n(R/R_0\cup R_\infty,d)}^\sim  
= I_n(R/R_0 \cup R_\infty,d)^\circ/S,$$
denoted by a superscripted tilde,
is determined by the (stack) quotient. The moduli space is 
empty unless $n>d$.
The rubber theory of $R$ is defined by integration against the
rubber virtual class,
 $$[{I_n(R/R_0\cup R_\infty,d)}^\sim ]^{vir}.$$ 
All of the above rubber constructions are $T$-equivariant.

The rubber moduli space $I_n(R/R_0\cup R_\infty, 0)^\sim$ carries
cotangent lines at the dynamical points $0$ and $\infty$ of $\Pp$. Let
$$\psi_0,  \psi_\infty \in A^1_T({I_n(R/R_0\cup R_\infty,d)}^\sim, {\mathbb Q})$$
denote the associated cotangent line classes.
Rubber integrals with relative conditions $\mu$ over $0$ and $\nu$
over $\infty$ are denoted by the bracket
\begin{equation} \label{rrc}
\left\langle \mu \ | \ \psi_0^a \psi_\infty^{b} \ |\ \nu \right\rangle_{n,d}^{
\sim}.
\end{equation}

Cotangent line classes in Donaldson-Thomas theory of $N$ can be
removed via topological recursion relations. For the
 relative theory  of $N/\cup_{i=1}^s N_{p_i}$, the topological recursion 
relation is:
\begin{multline*} 
\left\langle \left. \prod_i \sigma_{k_i}(F)  \
 \prod_{j=1}^s \psi_{p_j}^{a_j}  \ \right| \ \nu^1,\ldots, \nu^s
 \right\rangle_{n,d}^N =\\
\sum_{|\eta|=d}\ 
 \sum_{{n_1+n_2=n+d}} \left
\langle   \left.    \prod_i \sigma_{k_i}(F)\
\prod_{j\neq 1} \psi_{p_j}^{a_j}  \ \right| \ \eta, \nu^2, \ldots
\nu^s \right\rangle^N_{n_1,d} \\
\cdot \bigtriangleup_d(\eta,\eta)
\left. \left\langle \eta \ \right| \ \psi_\infty^{a_1-1} \ |\ \nu^1 \right\rangle_{n_2,d}^{\sim}.
\end{multline*}
The proof follows from the degeneration formula of Donaldson-Thomas theory
applied to the boundary expression for $\psi_{p_1}$ on the
Artin stack of target destabilizations.
The relative conditions away from $p_1$ and the descendent
insertions are bystanders in the topological recursion relation.

\subsection{Rubber calculus}
The rubber integrals \eqref{rrc} are determined via a Donaldson-Thomas
rubber calculus.
The technique, following Gromov-Witten theory \cite{mp,opvir}, involves rigidification
and topological recursion.

The universal target over the rubber moduli space is no longer a product. 
Let 
$$\pi:\mathcal{R} \rarr {I_n(R/R_0\cup R_\infty,d)}^\sim$$
denote the universal target.
The space $\mathcal{R}$ can be viewed as a moduli space of
rubber ideal sheaves {\em together} with a point $r$ of the target rubber.
The point $r$ is {\em not} permitted to lie on the relative divisors $R_0$ and
$R_\infty$. The stability condition is given by finiteness of the
associated automorphism group.
The virtual class of ${\mathcal R}$ is obtained via $\pi$-flat pull-back,
$$[{\mathcal R}]^{vir} = \pi^* \Big( [{I_n(R/R_0\cup R_\infty,d)}^\sim]^{vir}\Big).$$
As before, let
$${\mathfrak{J}} \rarr {\mathcal R}$$
denote the universal ideal sheaf on ${\mathcal R}$.

The target point $r$ together with $R_0$ and $R_\infty$
specifies 3 distinct points of the destabilized $\Pp$ over which the rubber is fibered.
By viewing the target point as $1\in \Pp$,
we obtain a rigidification map, 
$$\phi: \mathcal{R} \rarr I_n(N/N_0\cup N_\infty,d),$$
where $N=\cO_\Pp \oplus \cO_\Pp$ is the trivial bundle with 
level $(0,0)$ splitting over $\Pp$. By a comparison of deformation theories,
\begin{equation}\label{zek}
[{\mathcal R}]^{vir} = \phi^* \Big( [{I_n(N/N_0\cup N_\infty,d)}]^{vir}\Big).
\end{equation}

Rubber calculus transfers rubber integrals \eqref{rrc} to descendent integrals
on rigid (non-rubber) targets via the maps $\pi$ and $\phi$.
To start,
$$
(d-n) [I_n(R/R_0\cup R_\infty,d)]^{vir} =
\pi_*\Big( \text{ch}_3(\mathfrak{J}) \cap [{\mathcal R}]^{vir}\Big)$$
by a $\pi$-fiberwise calculation.
Since $n>d$, for nonempty rubber moduli spaces, $d-n$ is negative.
By the push-pull formula,
\begin{equation}\label{zex}
(d-n) 
\left\langle \mu \ | \ \psi_0^a \psi_\infty^{b} \ |\ \nu \right
\rangle_{n,d}^{\sim} =
\left\langle \mu \ |\ \text{ch}_3(\mathfrak{J}) \ \pi^*(\psi_0^a \psi_\infty^b) \ | \nu \right
\rangle_{n,d}^{{\mathcal R}\sim}.
\end{equation}
Next,
we compare the cotangent lines $\pi^*(\psi_0)$ and $\phi^*(\psi_0)$ on
${\mathcal R}$. A standard argument yields:
$$\pi^*(\psi_0)=\phi^*(\psi_0) - \phi^*(D_0),$$
where
 $$D_0 \subset I_n(N/N_0\cup N_\infty,d)$$
is the virtual boundary divisor
for which the rubber over $\infty$ carries Euler characteristic $n$.
Similarly,
 $$\pi^*(\psi_\infty)=\phi^*(\psi_\infty) - \phi^*(D_\infty).$$
We will apply the cotangent line comparisons to the right
side of \eqref{zex}.

Consider the Hilbert scheme of points ${\text{Hilb}}(R_0,d)$ of the
relative divisor.
The boundary condition $\mu$ corresponds to a Nakajima basis element
of $A^*_T({\text{Hilb}}(R_0,d), {\mathbb Q})$.
Let
${\mathfrak J}_0$ be the universal ideal sheaf on 
$$\text{Hilb}(R_0,d) \times R_0,$$ and let
$$\sigma_1=\pi_*\Big( {\text{ch}}_3({\mathfrak J}_0)\Big) \in A^1_T({\text{Hilb}}(R_0,d), 
{\mathbb Q}).$$
The class $\sigma_1$ in $A^1_T({\text{Hilb}}(R_\infty,d)$ is defined
in the same way.

The cotangent line comparisons and equation \eqref{zex} together yield
the following result:
\begin{multline}\label{dx}
(d-n) 
\left\langle \mu \ | \ \psi_0^a \psi_\infty^{b} \ 
|\ \nu \right\rangle_{n,d}^{\sim}  = \left\langle \mu \ |\ \sigma_1(F) \ \psi_0^a \psi_\infty^b \ |\ \nu \right
\rangle_{n,d}^{N} \\  \ \  \ \ \ \ \  \ \ \ \ \ \ \ \ \ \ \ \ \ \ \ \ \ \ 
 - \left\langle \sigma_1\cdot \mu \ | 
\ \psi_0^{a-1} \psi_\infty^{b} \ |\ \nu \right\rangle_{n,d}^{\sim} 
\\  - \left\langle \mu \ | 
\ \psi_0^{a} \psi_\infty^{b-1} \ |\ \sigma_1\cdot \nu \right\rangle_{n,d}^{\sim}. 
\end{multline}
The rubber integrals \eqref{rrc} are expressed in terms of the relative Donaldson-Thomas
theory of $N$ with descendent insertions $\sigma_k([F])$ by 
repeated applications of equation \eqref{dx} and the topological recursion
relations.

The following two vanishing statements are a consequence of the rubber calculus and
Lemma \ref{hx}.

\begin{Lemma} For $n>d(1-g)$,
\begin{equation*}
\left 
\langle \left.   \prod_{i} \sigma_{k_i}(F)\   \prod_{j=1}^s \psi_{p_j}^{a_j}\ \right| \ 
\nu^1,\ldots, \nu^s \right \rangle^{\cO_C\oplus \cO_C \pm}_{n,d} =0.
\end{equation*}
\end{Lemma}

\begin{Lemma} For $n>d$, \label{rubvan}
\begin{equation*}
\left. \Big 
\langle \mu \ \right| \ \psi_0^a \psi_\infty^b \ \left| \ \nu \Big \rangle\right.^{\sim\pm}_{n,d} =0.
\end{equation*}
\end{Lemma}

\subsection{Parallels}
The anti-diagonal vanishing in level $(0,0)$ holds for all vertices in the
triangle of equivalences of Section \ref{heq}.

The Gromov-Witten vanishing follows easily from Mumford's Hodge bundle relation,
$$c({\mathbb E}) \cdot c({\mathbb E^\vee}) = 1,$$
on the moduli space of curves $\overline{M}_g$, see \cite{bp}.
The Hilbert scheme vanishing is obtained from the existence of a 
modified virtual class in the hyperk\"ahler setting. The proof requires
a restriction of the obstruction theory of maps to $\text{Hilb}(\com^2,d)$, 
see \cite{qchs}.
The Donaldson-Thomas vanishing proven here could also be pursued
via a construction of a modified virtual class. Instead, our derivation 
proceeds formally from the equivariant vertex calculations of \cite{mnop1,mnop2}
using localization and degeneration.


\section{Additivity}
\subsection{Summary}
The level $(0,0)$ Donaldson-Thomas theory of $\Pp$ 
relative to $0,\infty\in \Pp$
will play a crucial role in the study of local curves.
The main results of the Section are vanishing
 and additivity properties for the invariants
\begin{equation}\label{cvgg}
\left \langle \mu \ \Big| \ \sigma_1(F) \ \Big| 
\  \nu\right\rangle.
\end{equation}
We follow here the bracket conventions of Section \ref{bcon}.

The vanishing of the invariants \eqref{cvgg} in most cases is
established by the following three results. 
\begin{Lemma} \label{xwq}
If $|\ell(\mu)-\ell(\nu)| >1$, then
$\left \langle \mu \ \Big| \ \sigma_1(F) \ \Big| 
\  \nu\right\rangle$ vanishes.
\end{Lemma}
\begin{Lemma} \label{xwqq}
If $|\ell(\mu)-\ell(\nu)| =1$, then
$\left \langle \mu \ \Big| \ \sigma_1(F) \ \Big| 
\  \nu\right\rangle_{n,d}$ vanishes
for $n>d$. 
\end{Lemma}

\begin{Lemma} \label{xwqqz}
If $\ell(\mu)=\ell(\nu)$, then
$\left \langle \mu \ \Big| \ \sigma_1(F) \ \Big| 
\  \nu\right\rangle$ vanishes
unless $\mu=\nu$. 
\end{Lemma}
In the diagonal case,
we will obtain the form 
\begin{equation}\label{kkwp}
\left \langle \mu \ \Big| \ \sigma_1(F) \ \Big| 
\  \mu\right\rangle_{n}= {\gamma_{\mu,n}}{(t_1t_2)^{-\ell(\mu)}}(t_1+t_2)
\end{equation}
for $\gamma_{\mu,n}
\in {\mathbb Q}$
and prove
a crucial additivity property parallel
to Equation (25) of \cite{qchs}.
\begin{Proposition} \label{sdds} 
An
additivity relation holds:
\begin{multline*}
  \label{addm} 
\frac{\left \langle \mu \ \Big|\ \sigma_1(F)\ \Big| \  \mu \ 
\right \rangle}
{\left \langle \mu \ \Big| \  \mu \right \rangle_{|\mu|,|\mu|}} =
\sum_i 
q^{|\mu|-\mu_i} \frac{\left \langle \mu_i \ \Big| \
 \sigma_1(F)\  \Big| \ \mu_i \right \rangle}
{\left \langle \mu_i \Big|\  \mu_i \right\rangle_{\mu_i,\mu_i}}\\ -
(\ell(\mu)-1) (t_1+t_2) q^{|\mu|} \Phi(q) \,.
\end{multline*}
\end{Proposition}
\noindent The bracket $\Big\langle \mu \ \Big|\ \mu\Big 
\rangle_{|\mu|,|\mu|}$ is 
the intersection form 
$$\Big\langle \mu\ \Big|\ \mu\Big\rangle_{|\mu|,|\mu|}
= \frac{1}{     (t_1 t_2)^{\ell(\mu)}         }
\frac{(-1)^{d-\ell(\mu)}}
{\mathfrak{z}(\mu)}$$
in the Nakajima basis. 
The function $\Phi(q)$ is defined by
$$\Phi(q) = q\frac{d}{dq} \log M(-q),$$
where $M(q)$ is the MacMahon series.

\subsection{Proofs of Lemmas \ref{xwq} and \ref{xwqq}}\label{ddaa}
Let $N=\cO_\Pp \oplus \cO_\Pp$ be the trivial bundle with
splitting of level $(0,0)$. Let $\overline{N}$ denote the
compactification of $N$ defined in Section \ref{thmx}.
The proofs  are obtained
from dimensional analysis for integrals in the Donaldson-Thomas
theory of $\overline{N}$.

By linearity, we can trade the invariants with unweighted
relative conditions for weighted relative conditions:
\begin{eqnarray*}\left \langle \mu \ \Big|
 \ \sigma_1(F) \ \Big| \  \nu\right\rangle & = &
\frac{1}{(t_1t_2)^{\ell(\mu)}}\left \langle \mu([0]) \ \Big| \ \sigma_1(F) 
\ \Big| \  \nu\right\rangle\\
& = &
\frac{1}{(t_1t_2)^{\ell(\nu)}}\left \langle \mu \ \Big| \ \sigma_1(F)
 \ \Big| \  \nu([0])\right\rangle
\end{eqnarray*}
By $T$-localization on $\overline{N}$ and the vanishing of Lemma \ref{feww}, 
we find
\begin{eqnarray*}
\left\langle\ \mu([0])\ \Big|\  \sigma_1(F)\ \Big| \  \nu\right
\rangle^{\overline{N}}_{n} & = & 
\left\langle\ \mu([0])\ \Big|\  \sigma_1(F)\ \Big| 
\  \nu\right\rangle_{n} \\
& &+ q^d\left\langle \mu([0])\ \Big|\  \nu\right\rangle
_{d,d} \cdot 
\left \langle \emptyset \ \Big|  \sigma_1(F) \Big| \ \emptyset \right
\rangle_{n-d}\Big|_{t_1-t_2,-t_2}\\
& &+ q^d\left\langle \mu([0])\ \Big|\  \nu\right\rangle
_{d,d} \cdot 
\left \langle \emptyset \ \Big| \sigma_1(F) \Big| \ \emptyset \right
\rangle_{n-d}\Big|_{t_2-t_1,-t_1}
\end{eqnarray*}
If $1+\ell(\mu)-\ell(\nu)<0$, then  
$$
\left\langle\ \mu([0])\ \Big|\  \sigma_1(F)\ \Big| \  \nu\right
\rangle^{\overline{N}}_{n}=0$$ by dimension
considerations since $\overline{N}$ is proper.
If $\ell(\mu)-\ell(\nu)<0$, then
$$\left\langle\ \mu([0])\ \Big|\  \nu\right\rangle_{d,d}=0.$$
Therefore, $1+\ell(\mu)-\ell(\nu)<0$ implies
$$
\left\langle\ \mu([0])\ \Big|\  \sigma_1(F)\ \Big| \  \nu\right
\rangle_{n}=0.$$
Similarly, if $1-\ell(\mu)+\ell(\nu)<0$, then
$$
\left\langle\ \mu\ \Big|\  \sigma_1(F)\ \Big| \  \nu([0])\right
\rangle_{n}=0.$$
We conclude
$\left \langle \mu \ \Big| \ \sigma_1(F) \ \Big| \  \nu\right\rangle$
vanishes unless
$|\ell(\mu)-\ell(\nu)|\leq 1$. Lemma \ref{xwq} is proven.\qed

If the equality $|\ell(\mu)-\ell(\nu)|=1$ holds, then the argument
yields a finer result. 
Either
$$
\left\langle\ \mu([0])\ \Big|\  \sigma_1(F)\ \Big| \  \nu\right
\rangle_{n,d} = 
\left\langle\ \mu([0])\ \Big|\  \sigma_1(F)\ \Big| \  \nu\right
\rangle^{\overline{N}}_{n,d}\ \in \mathbb{Q}
$$
or
$$
\left\langle\ \mu\ \Big|\  \sigma_1(F)\ \Big| \  \nu([0])\right
\rangle_{n,d} = 
\left\langle\ \mu\ \Big|\  \sigma_1(F)\ \Big| \  \nu([0])\right
\rangle^{\overline{N}}_{n,d}\in \mathbb{Q}.
$$
If $n>d$, the factor $(t_1+t_2)$ must divide the invariant
\begin{equation}
\label{fert}
\left \langle \mu \ \Big| \ \sigma_1(F) \ \Big| \  \nu\right\rangle_{n,d}
\end{equation}
by Lemma \ref{hx}. 
Thus, if $|\ell(\mu)-\ell(\nu)|=1$ and $n>d$, 
the integral \eqref{fert} must vanish.
Lemma \ref{xwqq} is proven.\qed

\subsection{Degree 0 calculation}
We calculate the integral
$\left \langle \emptyset \ \Big|\ \sigma_1(F)\ \Big| \ \emptyset \right
\rangle$.
By the degeneration formula and Theorem \ref{tm1},
\begin{eqnarray*}
\left \langle \emptyset \ \Big|\ \sigma_1(F)\ \Big| \ \emptyset \right
\rangle & = & \frac{\Big\langle \sigma_1(F) \Big\rangle_0}
{\left \langle \emptyset \ \Big| \right
\rangle \cdot \left \langle  \Big| \ \emptyset \right
\rangle}\\
& = &  M(-q)^{2\frac{(t_1+t_2)^2}{t_1t_2}} \Big\langle \sigma_1(F) 
\Big\rangle_0.
\end{eqnarray*}

To complete the calculation, we determine 
the series $\Big\langle \sigma_1(F) \Big\rangle_0$ via
localization with respect to the ${\mathbf T}$-action defined in
Section \ref{zloc}. Let the ${\mathbf T}$-equivariant lift of $F$
be specified by $[N_0]$. 
Then,
$$\Big\langle \sigma_1(F) \Big\rangle_0 =
\Big[\bW_{\sigma_1(s)}(\emptyset,\emptyset,\emptyset)|_{s,t_1,t_2}
\cdot \bW(\emptyset,\emptyset,\emptyset)|_{-s,t_1,t_2}\Big]_{s=0}.$$

\begin{Lemma} \label{cht} The vertex measure 
$\bW_{\sigma_1(s)}(\lambda,\emptyset,\emptyset)$ with a
descendent insertion is determined by
$$\bW_{\sigma_1(s)}(\lambda,\emptyset,\emptyset)|_{s,t_1,t_2}  = 
s\left(- q \frac{d}{dq} + c(\lambda;t_1,t_2)+ \frac{|\lambda|}{2}(t_1+t_2)  
\right) 
\bW (\lambda,\emptyset,\emptyset)|_{s,t_1,t_2}.$$
\end{Lemma}

Here, $c(\lambda; t_1,t_2)$ is the sum of the $(t_1,t_2)$-contents
of all squares in $\lambda$
$$
c(\lambda; t_1,t_2) = \sum_{(i,j)\in\lambda} \Big(
it_1 + jt_2 \Big) \,. 
$$
Viewing $\lambda$ as a Young diagram, 
the sum is over the {\em interior} corners of the squares --- the
corners closest to the origin.
For $|\lambda|\leq 1$, the total content $c(\lambda;t_1,t_2)$ vanishes.

\begin{proof}
Let ${\mathbf T}$ act of $\com^3$ with tangent weights $s,t_1,t_2$ at
the origin.
Let $\pi$ be a 3-dimensional partition with a single infinite leg of 
cross-section $\lambda$ in the direction of the tangent weight $s$.
 Let $I_\pi$ denote the corresponding $\mathbf{T}$-fixed
ideal. We first compute the restriction
$$\text{ch}_3({\mathfrak{I}})|_{[I_\pi]\times 0},$$
where $\mathfrak{I}$ is the universal ideal sheaf.

Let $s_1,s_2,s_3$ be the associated weights of the 
$\mathbf{T}$-action on the coordinates $x_1,x_2,x_3$ of $\com^3$,
$$s_1=-s, \ \ s_2=-t_1, \ \ s_3=-t_2.$$
Consider a graded free resolution of $I_\pi$ of
length $3$,
$$0 \rightarrow \bigoplus_{k} x^{c_k} A \rightarrow
\bigoplus_j x^{b_j} A  \rightarrow \bigoplus_{i} x^{a_i} A \rightarrow I_\pi 
\rightarrow 0,$$
where $A=\com[x_1,x_2,x_3]$ is the coordinate ring and 
$$
a_i, b_j, c_k \in \Z^3
$$
denote the degrees of the generators in each step.

Computing the Chern character via the resolution, we find
$$
\text{ch}_3({\mathfrak{I}})|_{I_\pi\times 0} = \frac{1}{3!} \left(
\sum_i (s, a_i)^3 
- \sum_j (s, b_j)^3 
+  \sum_k (s, c_k)^3 \right)\,.
$$ 
Here, $s=(s_1,s_2,s_3)$, and $(s, v)$ denotes the 
standard inner product of $s$ and $v$ in $\Z^3$. 

By calculating the trace of the $\mathbf{T}$-representation 
defined by $A/I_\pi$, we obtain a second relation:
\begin{multline}\label{resol_sum}
1-\sum_i e^{(s, a_i)} + \sum_j e^{(s, b_j)}
-\sum_k e^{(s, c_k)} =  \\
(1-e^{s_2})(1-e^{s_3}) 
\sum_{(i,j)\in\lambda} e^{is_2 + js_3} \, + \\
(1-e^{s_1})(1-e^{s_2}) (1-e^{s_3}) 
\sum_{p\in \pi'} e^{(s,p)}  
\end{multline}
where $\pi'$ denotes $\pi$ minus the infinite leg.
In particular, the renormalized volume $|\pi|$ is the number 
of squares in $\pi'$.

Extracting the cubic term in \eqref{resol_sum}, we find 
\begin{eqnarray*}
\text{ch}_3({\mathfrak{I}})|_{I_\pi\times 0} 
& = & - s_2 s_3 \left(c(\lambda; s_2,s_3)+\frac{|\lambda|}2\, 
(s_2+s_3)\right) + s_1 s_2 s_3 |\pi|\, \\
& = & t_1t_2 \left( -s |\pi| + c(\lambda;t_1,t_2) +\frac{|\lambda|}2\, 
(t_1+t_2)\right). 
\end{eqnarray*}
When applied to the computation of 
$\bW_{\sigma_1(s)}(\lambda,\emptyset,\emptyset)|_{s,t_1,t_2}$,
the argument of the descendent $\sigma_1(s)$ and the
equivariant push-forward together remove the
prefactor $t_1t_2$.
\end{proof}

By evaluation \eqref{fwq} of $\bW(\emptyset,\emptyset,\emptyset)$
and Lemma \ref{cht},
$$\bW_{\sigma_1(s)}(\emptyset,\emptyset,\emptyset) = 
\frac{(t_1+t_2)(t_1+s)(t_2+s)}{t_1t_2}\ \Phi(q)\
\bW(\emptyset,\emptyset,\emptyset)|_{s,t_1,t_2}.
$$
After multiplying all the factors, 
\begin{equation}\label{dggz}
\Big\langle \emptyset \ \Big| \ \sigma_1(F)\ \Big| \ 
\emptyset \Big\rangle = (t_1+t_2)\ \Phi(q).
\end{equation}

\subsection{Proofs of Lemma \ref{xwqqz} and Proposition \ref{sdds}}
The proof is via $T$-localization on $\overline{N}$.
We follow the notation of Section \ref{thmx} for the $T$-fixed
points of $\overline{N}_\infty$. Let
$$\fpa, \fpb\in \overline{N}_\infty$$ be the two new
$T$-fixed points with normal
$T$-weights $$t_1-t_2, -t_2, \ \text{and} \ t_2-t_1,-t_1$$
respectively.

For $\ell(\mu)=\ell(\nu)$,
consider the $T$-equivariant localization formula for
the integral
\begin{eqnarray*}
\left\langle\ \mu \ \Big|\  \sigma_1(F)\ \Big| \  \nu([0])\right
\rangle^{\overline{N}} & = & \quad
\left\langle\ \mu \ \Big|\  \sigma_1(F)\ \Big| 
\  \nu([0]) \right\rangle \\
& &+ q^d\left\langle\ \mu \ \Big|\  \nu([0])\right\rangle
_{d,d}\ \cdot \
\left \langle \emptyset \ \Big| \ \sigma_1(F) \ \Big| \ \emptyset \right
\rangle\Big|_{t_1-t_2,-t_2}\\
& &+ q^d\left\langle\ \mu\ \Big|\  \nu([0])\right\rangle
_{d,d}\ \cdot \
\left \langle \emptyset \ \Big| \ \sigma_1(F) \ \Big| \ \emptyset \right
\rangle\Big|_{t_2-t_1,-t_1}.
\end{eqnarray*}
After evaluation and rearrangement,
\begin{multline}\label{kdd}
\left\langle\ \mu \ \Big|\  \sigma_1(F)\ \Big| 
\  \nu([0]) \right\rangle =
\left\langle\ \mu \ \Big|\  \sigma_1(F)\ \Big| \  \nu[(0)]\right
\rangle^{\overline{N}}\\
+      \frac{(-1)^{d-\ell(\mu)} \delta_{\mu,\nu}}
{\mathfrak{z}(\mu)}        (t_1+t_2)  q^d\Phi(q) 
\end{multline}
We find 
$\left\langle\ \mu \ \Big|\  \sigma_1(F)\ \Big| 
\  \nu[(0)] \right\rangle$
is a linear function of
$t_1$ and $t_2$. By symmetry, the integral must be a function
of $t_1+t_2$. 

We draw two conclusions. First, in the diagonal case $\mu=\nu$,
the form \eqref{kkwp} is proven. Second,
$\left\langle\ \mu \ \Big|\  \sigma_1(F)\ \Big| 
\  \nu[(0)] \right\rangle$
is {\em determined} by restriction to $t_2=0$.

If $\ell(\nu)=1$, the Lemma and Proposition are
empty. Let $\ell(\nu)\geq 2$. We proceed by induction on
$\ell(\nu)$.

The strategy of the proof is to compare two $T$-equivariant
integrals:
$$
m(\nu_1)\left\langle\ \mu \ \Big|\  \sigma_1(F)\ \Big| \  \nu[(0)] \right
\rangle^{\overline{N}}_{n}$$
and
$$
\left\langle\ \mu \ \Big|\  \sigma_1(F)\ \Big| \ \nu_1([\fpa]),
\nu_2[(0)], \ldots,
\nu_{\ell(\nu)}([0])\right
\rangle^{\overline{N}}_{n},$$
where $m(\nu_1)$ is the multiplicity of the part $\nu_1$ in $\nu$.
Both integrals 
evaluate to linear functions of $t_1$ and $t_2$. The answers
are {\em not} equal --- the integrands are different $T$-equivariant
classes. Consider the 1-dimensional subtorus
$$T_1\subset T$$
determined by the first factor. The points 
$$0,\fpa \in \overline{N}_\infty$$
define {\em equivalent} $T_1$-equivariant classes on $\overline{N}$.
Hence, the two integrals above are equal after the restriction
$t_2=0$. 

Consider the $T$-equivariant localization formula for
the second integral. The formula immediately yields 0
unless there is a part $\mu_1$ of $\mu$ equal to $\nu_1$. Let
$\mu^*$ and $\nu^*$ denote the subpartitions obtained
by removing the first parts of $\mu$ and $\nu$. Then, 
\begin{multline*}
\left\langle\ \mu \ \Big|\  \sigma_1(F)\ \Big| \  \nu_1([\fpa]),\nu^*([0])
\right
\rangle^{\overline{N}} =\\ \hspace{-10pt}
q^{\mu_1} \frac{(-1)^{\mu_1-1}}{\mu_1} 
\left\langle\ \mu^* \ \Big|\  \sigma_1(F)\ \Big| 
\  \nu^*([0]) \right\rangle \\ \hspace{+105pt}
+ q^{|\mu^*|} \frac{(-1)^{|\mu^*|-\ell(\mu^*)} \delta_{\mu^*,\nu^*}}
{\mathfrak{z}(\mu^*)}
{\left\langle\ \mu_1 \ \Big|\ \sigma_1(F) \ \Big| \  \nu_1([0])\right\rangle}
\Big|_{t_1-t_2,-t_2}\\ \hspace{-15pt}
+ q^{d}
\frac{(-1)^{\mu_1-1}}{\mu_1} \frac{(-1)^{|\mu^*|-\ell(\mu^*)} 
\delta_{\mu^*,\nu^*}}
{\mathfrak{z}(\mu^*)} (-2t_1+t_2)  \Phi(q)
\end{multline*}
Comparing the $t_2=0$ restriction with \eqref{kdd} and using
induction, we find 
$$\left\langle\ \mu \ \Big|\  \sigma_1(F)\ \Big| 
\  \nu[(0)] \right\rangle=0$$
unless $\mu=\nu$ proving Lemma \ref{xwqqz}.

We now assume $\mu=\nu$. Combining all the equations yields the
following additivity relation:
\begin{multline*}
\frac{\left\langle\ \mu \ \Big|\  \sigma_1(F)\ \Big| 
\  \mu \right\rangle}{\left\langle\ \mu \ \Big| \ \mu \ \right 
\rangle_{|\mu|,|\mu|}}  =  \quad 
q^{\mu_1} 
\frac{\left\langle\ \mu^* \ \Big|\  \sigma_1(F)\ \Big| 
\  \mu^* \right\rangle}
{\left\langle\ \mu \ \Big| \ \mu \ \right 
\rangle_{|\mu^*|,|\mu^*|}}
\\
+ q^{|\mu^*|} 
\frac{
\left\langle\ \mu_1 \ \Big|\ \sigma_1(F) \ \Big| \  \mu_1\right\rangle}
{\left\langle\ \mu_1 \ \Big| \ \mu_1 \ \right 
\rangle_{\mu_1,\mu_1}}
- (t_1+t_2) q^{d} \Phi(q).
\end{multline*}
The induction step is complete and Proposition \ref{sdds} is proven.
\qed

\section{Degree 1}\label{ddd111}
\subsection{Vertex measure}\label{ttmm}
The following formula for the 
degree 1 vertex measure was stated (without proof) in
\cite{mnop2}. 
\begin{Proposition} \label{xccv} At a 3-fold fixed point with tangent weights
$s_1$, $s_2$, and $s_3$, the degree 1 vertex measure is
$$\bW(1,\emptyset,\emptyset)=   (1+q)^{\frac{s_2+s_3}{s_1}} 
M(-q)^{- \frac{(s_1+s_2)(s_1+s_3)(s_2+s_3)}{s_1s_2s_3}},$$
where the unique degree 1 leg extends in the $s_1$ direction.
\end{Proposition}

As an easy application of Proposition \ref{xccv}, we calculate the
series
$$
\left \langle\   \right
\rangle_1 =
\sum_{n\geq 1} q^n \left \langle\   \right \rangle_{n,1}$$
via ${\mathbf T}$-localization:
\begin{eqnarray*}
\left \langle\   \right
\rangle_1& = & \Big[\bW(1,\emptyset,\emptyset)|_{s,t_1,t_2}
\cdot \frac{q}{t_1t_2}\cdot  \bW(1,\emptyset,\emptyset)|_{-s,t_1,t_2}\Big]_{s=0}\\
& = &\frac{q}{t_1t_2}  M(-q)^{-2\frac{(t_1+t_2)^2}{t_1t_2}}.
\end{eqnarray*}
Empty Donaldson-Thomas brackets $\left\langle \  \right\rangle$ denote the
integrand $1$,
$$\left \langle\   \right
\rangle_1 = \left\langle  1 \right \rangle_1.$$

\subsection{Proof of Proposition \ref{xccv}}
The Proposition is proven by using two geometric constraints for 
the vertex measure $\bW(1,0,0)$.

The effective curve classes of $\Pp \times \Pt$ are generated by
$\beta_1$ and $\beta_2$ where
\begin{eqnarray*}
\beta_1 & = & [\Pp \times 0], \\
\beta_2 & = & [0 \times \Pp].
\end{eqnarray*}
We will calculate Donaldson-Thomas series associated to
these two classes.

Consider first $\beta_1$. The virtual dimension of $I_n(\Pp \times
\Pt,\beta_1)$ is 2. There is a Hilbert-Chow morphism
$$\epsilon:I_n(\Pp \times \Pt,\beta_1) \rarr \Pt$$
which specifies the location of the line over the second factor
of $\Pp \times \Pt$. We will compute the Donaldson-Thomas
series
\begin{equation}\label{bbggr}
\frac{\Big \langle   \epsilon^*(P) \Big \rangle^{\Pp\times
\Pt}_{\beta_1}}
{\Big \langle   \Big \rangle^{\Pp\times
\Pt}_{0}}
 =\frac{
\sum_{n\geq 1} q^n 
\Big \langle  \epsilon^*(P) \Big \rangle^{\Pp \times \Pt}
_{n,\beta_1}}
{
\sum_{n\geq 0} q^n 
\Big \langle  \  \Big \rangle^{\Pp \times \Pt}
_{n,0}}
,
\end{equation}
where $P$ is  the class of a point in $\Pt$.
The integrand in the numerator of
\eqref{bbggr} has dimension 2, so the 
integrals are well-defined.

We will calculate \eqref{bbggr} by equivariant localization.
Let the $1$-dimensional torus
$S$ act on $\Pp$ with
fixed points $0,\infty\in \Pp$ and tangent weights $s,-s$.
Let the $2$-dimensional torus $T$ act on $\Pt$ with fixed points
$p,p',p''$. Let $t_1,t_2$ be the tangent $T$-weights at $p$.
Localization of \eqref{bbggr} by the action of the
3-dimensional torus,
$$\mathbf{T}=S \times T,$$
yields
$$ 
\frac{\Big \langle   \epsilon^*(P) \Big \rangle^{\Pp\times
\Pt}_{\beta_1}}
{\Big \langle \ \Big \rangle^{\Pp\times
\Pt}_{0}}
=
t_1t_2\frac{\bW(1,\emptyset,\emptyset)}{\bW(\emptyset,\emptyset,\emptyset)}
\Big|_{s,t_1,t_2}\cdot \frac{q}{t_1t_2} \cdot  
\frac{\bW(1,\emptyset,\emptyset)}{\bW(\emptyset,\emptyset,\emptyset)}
\Big|_{-s,t_1,t_2}.$$
Here,
the $\mathbf{T}$-equivariant lift of $P$ is specified by $P=[p]$.
We conclude
\begin{equation}\label{xxon}
{\bW'}(1,\emptyset,\emptyset)|_{s,t_1,t_2}
{\bW'}(1,\emptyset,\emptyset)|_{-s,t_1,t_2}\in \mathbb{Q}[[q]],
\end{equation}
where
$$
{\bW'}(1,\emptyset,\emptyset)=
\frac{\bW(1,\emptyset,\emptyset)}{\bW(\emptyset,\emptyset,\emptyset)}.
$$

Next, consider the class $\beta_2$. The virtual dimension of $I_n(\Pt \times
\Pp,\beta_1)$ is 3. There is a Hilbert-Chow morphism
$$\epsilon:I_n(\Pp \times \Pt,\beta_2) \rarr \Pp \times 
(\mathbf{P}^{2})^\vee$$
which specifies the line component.
We will compute the Donaldson-Thomas
series
\begin{equation}\label{bbggrr}
\frac{\Big \langle   \epsilon^*(Q) \Big \rangle^{\Pp\times
\Pt}_{\beta_2}}
{\Big \langle   \  \Big \rangle^{\Pp\times
\Pt}_{0}}
 =\frac{
\sum_{n\geq 1} q^n 
\Big \langle  \epsilon^*(Q) \Big \rangle^{\Pp \times \Pt}
_{n,\beta_2}}
{
\sum_{n\geq 0} q^n 
\Big \langle  \  \Big \rangle^{\Pp \times \Pt}
_{n,0}}
,
\end{equation}
where $Q$ is  the class of a point in $\Pp\times (\mathbf{P}^{2})^\vee$.
The integrand in the numerator of
\eqref{bbggr} has dimension 3, so the 
integrals are well-defined.

The ${\mathbf T}$-equivariant 
localization formula yields the following evaluation of the
series \eqref{bbggrr}:
$$ 
\frac{\Big \langle   \epsilon^*(Q) \Big \rangle^{\Pp\times
\Pt}_{\beta_2}}
{\Big \langle  \  \Big \rangle^{\Pp\times
\Pt}_{0}}
=
st_2(t_2-t_1)\frac{\bW(1,\emptyset,\emptyset)}{\bW(\emptyset,\emptyset,\emptyset)}
\Big|_{t_1,s,t_2}\cdot \frac{q}{st_2(t_2-t_1)} \cdot  
\frac{\bW(1,\emptyset,\emptyset)}{\bW(\emptyset,\emptyset,\emptyset)}
\Big|_{-t_1,s,t_2-t_1}.$$
The $\mathbf{T}$-equivariant lift of $Q$ is specified by the line
over $0\in \Pp$ connecting $T$-fixed points $p,p'\in \Pt$ where
the tangent $T$-weights at $p'$ are $-t_1, t_2-t_1$.
We conclude
\begin{equation*}
{\bW'}(1,\emptyset,\emptyset)|_{t_1,s,t_2}
{\bW'}(1,\emptyset,\emptyset)|_{-t_1,s,t_2-t_1}\in \mathbb{Q}[[q]].
\end{equation*}
After renaming the variables, we obtain
\begin{equation}\label{xxtwo}
{\bW'}(1,\emptyset,\emptyset)|_{s,t_1,t_2}
{\bW'}(1,\emptyset,\emptyset)|_{-s,t_1,t_2-s}\in \mathbb{Q}[[q]],
\end{equation}
and, by symmetry,
\begin{equation}\label{xxtwop}
{\bW'}(1,\emptyset,\emptyset)|_{s,t_1,t_2}
{\bW'}(1,\emptyset,\emptyset)|_{-s,t_1-s,t_2}\in \mathbb{Q}[[q]],
\end{equation}

By definition of the equivariant vertex measure,
$$\bW(1,\emptyset,\emptyset)|_{s,t_1,t_2} \in {\mathbb Q}(s,t_1,t_2)[[q]].$$
The $q^0$ coefficient is 1.
By Lemma \ref{xtyy}, the coefficient of $q^n$ is divisible by
$t_1+t_2$ for $n>0$.

By repeated applications of the 
logarithms of equations \eqref{xxon}, \eqref{xxtwo}, and
\eqref{xxtwop}, 
we find
$$
\log \Big({\bW'}(1,\emptyset,\emptyset)\Big)\Big|_{s,t_1,t_2}  +
\log \Big(
{\bW'}(1,\emptyset,\emptyset)\Big)
\Big|_{-s,t_1-is,t_2-js}\in \mathbb{Q}[[q]]
$$
for all non-negative integers $i$ and $j$.
The coefficients
$$
\log \Big({\bW'}(1,\emptyset,\emptyset)\Big)\Big|_{s,t_1,t_2} =
\sum_{n\geq 1} f_n(s,t_1,t_2) q^n$$
must therefore satisfy 
$$f_n(s,t_1,t_2)+ f_n(-s,t_1-xs,t_2-ys)=g_n(x,y)$$
for variables $x$ and $y$.
Differentiation with respect to $x$ yields
$$-s\frac{\partial f_n}{\partial t_1}(-s,t_1,t_2)= 
\frac{\partial g_n}{\partial x}(0,0).$$
Similarly
$$-s\frac{\partial f_n}{\partial t_2} (-s,t_1,t_2)= 
\frac{\partial g_n}{\partial y}(0,0).$$
Hence, by integration and symmetry,
$$f_n(s,t_1,t_2)= \gamma_n\frac{t_1+t_2}{s},$$
where $\gamma_n\in {\mathbb Q}$. Since $f_n$ {\em must} be
divisible by $t_1+t_2$, the constant of integration 
vanishes.

After specializing to Calabi-Yau weights $s+t_1+t_2=0$, we find
$$f_n= -\gamma_n.$$
However, after Calabi-Yau specialization, the
equivariant vertex measure takes the simple form
$$\bw(\pi)= (-1)^{|\pi|}$$
for 3-dimensional partitions $\pi$, see \cite{mnop1}.

\begin{Lemma}\label{dfix}
 $\log\Big({\bW'}(1,\emptyset,\emptyset)\Big)
\Big|_{s+t_1+t_2=0} =- \log(1+q)$.
\end{Lemma}

The Lemma
is a special case of the vertex evaluation 
required for the calculation of the level $(-1,0)$ cap \cite{OR},
see Section \ref{ver12} for a detailed discussion.

Lemma \ref{dfix} determines the constants $\gamma_n$ and completes
the proof of Proposition \ref{xccv}. \qed

\subsection{Descendent calculation}
We calculate the series
$\left \langle (1) \ \Big|\ \sigma_1(F)\ \Big| \ (1) \right
\rangle$.
By the degeneration formula,
\begin{eqnarray*}
\left \langle (1) \ \Big|\ \sigma_1(F)\ \Big| \ (1) \right
\rangle & = & \frac{q^2}{(t_1t_2)^2}
\frac{\Big\langle \sigma_1(F) \Big\rangle_1}
{\left \langle (1) \ \Big| \right
\rangle \cdot \left \langle  \Big| \ (1) \right
\rangle}\\
& = &  \frac{q}{t_1t_2}
\frac{\Big\langle \sigma_1(F) \Big\rangle_1}
{ \left \langle\   \right
\rangle_1}\\
& = & M(-q)^{2\frac{(t_1+t_2)^2}{t_1t_2}} \Big\langle \sigma_1(F) 
\Big\rangle_1.
\end{eqnarray*}

To complete the calculation, we determine 
the series $\Big\langle \sigma_1(F) \Big\rangle_1$ via
localization with respect to the ${\mathbf T}$-action.
Let the ${\mathbf T}$-equivariant lift of $F$
be specified by $[N_0]$. Then,
$$\Big\langle \sigma_1(F) \Big\rangle_1 =
\Big[\bW_{\sigma_1(s)}(1,\emptyset,\emptyset)|_{s,t_1,t_2}
\cdot 
\frac{q}{t_1t_2}\cdot
 \bW(1,\emptyset,\emptyset)|_{-s,t_1,t_2}\Big]_{s=0}.$$
%
By Proposition \ref{xccv} and  
Lemma \ref{cht},
\begin{eqnarray*}
\bW_{\sigma_1(s)}(1,\emptyset,\emptyset)
 & = & 
(t_1+t_2)(\frac{-q}{1+q}+\frac{1}{2}) \bW(1,\emptyset,\emptyset)_{s,t_1,t_2}\\& & \quad +
\frac{(t_1+t_2)(t_1+s)(t_2+s)}{t_1t_2} \Phi(q)
\bW(1,\emptyset,\emptyset)|_{s,t_1,t_2}.
\end{eqnarray*}
After multiplying all the factors, 
\begin{equation}\label{deg1}
\frac{\Big\langle (1) \ \Big| \ \sigma_1(F)\ \Big| \ 
(1) \Big\rangle}{
\left \langle (1) \ \Big| \ (1) \right\rangle_{1,1}} = \frac{(t_1+t_2)}{2}
q
\frac{1-q}{1+q} + (t_1+t_2) q\Phi(q).
\end{equation}
\label{cvbbn}

\subsection{Cap}
We calculate the partition function $\mathsf{DT}(0|0,0)_\lambda$
corresponding to the cap in the TQFT formalism.
 
\begin{Lemma}\label{lcap}
The invariants of the level $(0,0)$ cap are given by
\[
{\mathsf{DT}} (0| 0,0)_{\lambda }=\begin{cases}
\frac{1}{d!(t_{1}t_{2})^{d}}& \text{if $\lambda= (1^{d})$}\\
0&\text{if $\lambda \neq (1^{d})$.}
\end{cases}
\]
\end{Lemma}

\begin{proof} 
If $\lambda \neq (1^d)$, then $\ell(\lambda)<d$. By equation
\eqref{van7},
\begin{equation*}
\Big \langle \ 
\Big| \ \lambda([0]) \Big \rangle=0
\end{equation*}
if $\ell(\lambda)<d$.
Since
\begin{eqnarray*}
{\mathsf{DT}} (0| 0,0)_{\lambda } & = & 
q^{-d} \frac{
\Big\langle  \ 
\Big| \ \lambda \Big \rangle}
{
\Big\langle  \ 
\Big| \ \emptyset \Big \rangle}
 \\
& = & 
\frac{q^{-d}}{
(t_1t_2)^{\ell(\lambda)}} 
\frac{
\Big \langle  \ 
\Big| \ \lambda([0]) \Big \rangle}
{
\Big\langle  \ 
\Big| \ \emptyset \Big \rangle}\ ,
\end{eqnarray*}
the cap ${\mathsf{DT}} (0| 0,0)_{\lambda}$
vanishes if $\lambda\neq (1^d)$.

If $\lambda=(1^d)$, then the Donaldson-Thomas invariant
for the compactified geometry
\begin{equation*}
\Big \langle \ 
\Big| \ \lambda([0]) \Big \rangle^{\overline{(0,0)}}
\end{equation*}
is a constant independent of the equivariant parameters $t_1$ and $t_2$.
By specializing the partition weights to different
$T$-fixed points and localizing, we obtain
\begin{equation*}
\frac{
\Big \langle \ 
\Big| \ \lambda([0]) \Big \rangle}
{\Big \langle \ 
\Big| \ \emptyset \Big \rangle} = \frac{1}{d!}
\left( \frac{
\Big \langle \ 
\Big| \ 1([0]) \Big \rangle}
{\Big \langle \ 
\Big| \ \emptyset \Big \rangle} \right)^d.
\end{equation*}
The degree 1 calculation,
$$\frac{
\Big \langle \ 
\Big| \ 1([0]) \Big \rangle}
{\Big \langle \ 
\Big| \ \emptyset \Big \rangle} =q,$$
completes the proof of the Lemma.
\end{proof}

\section{The operator $\MM_\sigma$}\label{s2} 

\subsection{Fock space formalism}
\label{fs}
Let the $2$-dimensional torus $T$ act on $\com^2$ by standard
diagonal scaling.
We review the Fock space description of the $T$-equivariant 
cohomology of the Hilbert scheme of points of $\C^2$, see \cite{Gron, Nak}. 

By definition, the Fock space $\cF$ 
is freely generated over $\Q$ by commuting 
creation operators $\alpha_{-k}$, $k\in\Z_{>0}$,
acting on the vacuum vector $\vac$. The annihilation 
operators $\alpha_{k}$, $k\in\Z_{>0}$, kill the vacuum 
$$
\alpha_k \cdot \vac =0,\quad k>0 \,,
$$
and satisfy the commutation relations
$$
\left[\alpha_k,\alpha_l\right] = k \, \delta_{k+l}\,. 
$$

A natural basis of $\cF$ is given by 
the vectors  
\begin{equation}
  \label{basis}
  \lv \mu \rang = \frac{1}{\zz(\mu)} \, \prod \alpha_{-\mu_i} \, \vac \,.
\end{equation}
indexed by partitions 
$\mu$. Here, $$\zz(\mu)=|\Aut(\mu)| \, \prod \mu_i$$ is the usual 
normalization factor.

The { Nakajima basis} defines a canonical isomorphism,
\begin{equation*}
\cF \otimes _{\mathbb Q} {\mathbb Q}[t_1,t_2]\stackrel{\sim}{=} 
\bigoplus_{d\geq 0} A_{T}^*(\Hilb(\com^2,d),{\mathbb Q}).
\end{equation*}
The Nakajima basis element corresponding to  
$\lv \mu \rang$  is
$$\frac{1}{\Pi_i \mu_i} [V_\mu]$$
where $[V_\mu]$ is (the cohomological dual of)
the class of the subvariety of $\Hilb(\C^2,|\mu|)$
with generic element given by a union of 
schemes of lengths $$\mu_1, \ldots, \mu_{\ell(\mu)}$$ supported
at $\ell(\mu)$ distinct points of $\C^2$. 
The vacuum vector $\vac$ corresponds to the unit in 
$A_T^*(\Hilb_0, {\mathbb Q}).$

The standard inner product on the $T$-cohomology  induces the following 
{\em nonstandard} inner product on Fock space after an extension of scalars:
\begin{equation}
  \label{inner_prod}
 \Big \langle \mu \ \Big|\  \nu \Big\rangle_{\mathcal F} 
= \frac{1}{(t_1 t_2)^{\ell(\mu)}} 
\frac{(-1)^{|\mu|-\ell(\mu)}} {\zz(\mu)}
\ {\delta_{\mu\nu}} \,. 
\end{equation}
With respect to the inner product, 
\begin{equation}
  \label{adjoint}
  \left(\alpha_{k}\right)^* = (-1)^{k-1} (t_1 t_2)^{\sgn(k)} \, 
\alpha_{-k}\,.
\end{equation}
If there is no ambiguity, the subscript $\mathcal{F}$ will be
omitted from the bracket \eqref{inner_prod}.

\subsection{The class $D$}\label{cbns}
Let $\cO/{\cI}$ be the rank $d$ tautological bundle
on $\text{Hilb}(\com^2,d)$, and let
$$D= c_1(\cO/{\cI}).$$
A straightforward calculation shows 
$$
D= - \lv 2,1^{d-2} \rang \,.
$$
The classical multiplication of $D$ on the Fock space 
${\mathcal F}$ is given by the following operator:
\begin{equation}
-(t_1+t_2) \sum_{k>0} \frac{k-1}{2}  \,
 \alpha_{-k} \, \alpha_k  
+ \frac12 \sum_{k,l>0} \label{hhj}
\Big[t_1 t_2 \, \alpha_{k+l} \, \alpha_{-k} \, \alpha_{-l} -
 \alpha_{-k-l}\,  \alpha_{k} \, \alpha_{l} \Big],
\end{equation}
see \cite{Lehn,LS,LQW,qchs}.

The first summand of \eqref{hhj}
 contains a term proportional to the energy operator, 
$$
|\cdot| = \sum_{k>0} \alpha_{-k} \, \alpha_k \,.
$$
The energy operator
acts diagonally
on Fock space with eigenvalue $|\mu|$ on
the vector $|\mu\rangle$. 

\subsection{Operators}
\label{t1}

The following operator 
on Fock space plays a central role in the paper: 
\begin{multline}
  \label{theM} 
\MM(q,t_1,t_2) = (t_1+t_2) \sum_{k>0} \frac{k}{2} \frac{(-q)^k+1}{(-q)^k-1} \,
 \alpha_{-k} \, \alpha_k  + \\
\frac12 \sum_{k,l>0} 
\Big[t_1 t_2 \, \alpha_{k+l} \, \alpha_{-k} \, \alpha_{-l} -
 \alpha_{-k-l}\,  \alpha_{k} \, \alpha_{l} \Big] \,.
\end{multline}
The $q$-dependence of $\MM$ is only in the first sum in \eqref{theM}
which acts diagonally in the basis \eqref{basis}. The two terms in the 
second sum in \eqref{theM} are known respectively as the splitting and joining terms. 
The operator $\MM$ is self-adjoint
\begin{equation}
  \label{Madj}
  \MM^* = \MM 
\end{equation}
with respect to \eqref{adjoint}.

Let the operator ${\MM}_\sigma$ on
Fock space be defined by matrix elements
\begin{eqnarray*}
\Big \langle \mu \ \Big| \ {\MM}_\sigma \ \Big| \ 
\nu \Big \rangle_{\mathcal F} 
& = &
q^{-d} \Big \langle \mu \ \Big| -\sigma_1(F)\ \Big|\ \nu \Big \rangle\\
& = & \sum_{n\geq d} q^{n-d} 
\Big \langle \mu \ \Big|  -\sigma_1(F)\ \Big|\ \nu \Big \rangle_{n,d}
\end{eqnarray*}
for partitions satisfying $$|\mu|=|\nu|=d.$$
The matrix elements of ${\MM}_\sigma$ are understood
to vanish unless $|\mu|=|\nu|$.

\begin{Proposition}\label{T12}$
    {\MM}_\sigma = \MM - (t_1+t_2) \Phi(q) \cdot \text{\em Id}.$
\end{Proposition}

We check here three initial compatibilities required for
Proposition \ref{T12}. Proposition \ref{T12} will be proven in 
 Section 
\ref{prt} by relating Donaldson-Thomas integrals to Gromov-Witten
invariants of the Hilbert scheme $\text{Hilb}(\com^2,d)$.

First, Proposition \ref{T12} requires the negative of 
the $q$-shifted  operator bracket,
$$-q^d \Big \langle \mu \ \Big| \ {\MM}- (t_1+t_2)\Phi(q) \ \Big|\ \nu \Big 
\rangle_{{\mathcal F}}\,,$$ 
to satisfy the additivity property of Proposition \ref{sdds}.
The required additivity is easily checked.

Second, Proposition \ref{T12} is valid for $q=0$.
Let $N$ be the 
trivial bundle with level $(0,0)$ splitting on $\Pp$.
There is an isomorphism of moduli spaces
\begin{equation}\label{xfg}
I_d(N/N_0\cup N_\infty,d) \stackrel{\sim}{=} \text{Hilb}(N_z,d)
\end{equation}
for any $z\in \Pp$. 
Under the isomorphism \eqref{xfg}, the Donaldson-Thomas
descendent class
$-\sigma_1([N_z])$ determines an element of  
$A^1_T(\text{Hilb}(N_z,d),{\mathbb Q})$. By a Riemann-Roch
calculation,
$$-\sigma_1(N_z) = D -\frac{t_1+t_2}{2}d.$$ 
Proposition \ref{T12} correctly equates the $q=0$
evaluation ${\MM}_{\sigma}(0)$ of  ${\MM}_\sigma$
with    
the classical multiplication \eqref{hhj} of $D$ shifted
by $-\frac{t_1+t_2}{2} d$.

Third, Proposition \ref{T12} is valid in degrees 0 and 1
 by  the descendent 
evaluation \eqref{dggz}, 
$$\Big \langle \emptyset \ \Big| \ {\MM}_\sigma \ 
\Big| \ \emptyset \Big \rangle_{\mathcal F} =
\Big \langle\emptyset \ \Big| \ -\sigma_1(F) \ \Big| \emptyset \Big \rangle
= -(t_1+t_2) \Phi(q),$$
and the descendent evaluation
\eqref{deg1},
\begin{eqnarray*}
\frac{\Big \langle (1) \ \Big| \ {\MM}_\sigma \ 
\Big| \ (1) \Big \rangle_{\mathcal F}}{ 
\Big \langle (1) \  
\Big| \ (1) \Big \rangle_{\mathcal F} }
& = & 
q^{-1}\frac{\Big\langle (1) \ \Big| \ -\sigma_1(F)\ \Big| \ 
(1) \Big\rangle}{
\left \langle (1) \ \Big| \ (1) \right\rangle_{1,1}}\\
& =& -\frac{(t_1+t_2)}{2}
\frac{1-q}{1+q} - (t_1+t_2) \Phi(q).
\end{eqnarray*}

\section{Proof of Proposition \ref{T12}}
\label{prt}

\subsection{Induction strategy}

Our proof of Proposition \ref{T12} closely follows the 
proof of Theorem 1 in \cite{qchs} though several differences
occur. Since the results and the geometry
differ, the closeness of the proof
is somewhat surprising.

We proceed by induction on the degree $d$. If $d=0$ or $1$, 
Proposition \ref{T12} has already been proven by descendent
calculations.
Let $d>1$. 

Next, we induct on the
Euler number $n$. In the minimal case $n=d$, Proposition \ref{T12} recovers
classical multiplication by $D$ on the Hilbert scheme of points.
Let $n>d$.

The induction step relies upon the  
addition formula of Proposition \ref{sdds}.
We will compute 
an invariant
$$
\Big \langle \gamma_1 \ \Big| -\sigma_1(F)\ \Big| \  \gamma_2 \Big \rangle
_{n}$$
for which the expansions of the classes
$$\gamma_1,\gamma_2 \in A^{2d}_T(\Hilb(\com^2,d),{\mathbb Q}),$$
in the Nakajima basis contain nontrivial multiples 
{\em not divisible by $(t_1+t_2)$}
of the class $|(d)\rangle$.
By the addition rules, if 
\begin{equation}
  \label{t123}
\Big  \langle \gamma_1 \ \Big| \ \MM_\sigma \ \Big| \ \gamma_2 \Big
\rangle_{{\mathcal F},n-d} = 
\Big \langle  \gamma_1 \ \Big| -\sigma_1(F)\ \Big| \  \gamma_2 
\Big \rangle_{n} \,,
\end{equation}
then Proposition \ref{T12} is proven for Euler number $n$. 

Both sides of 
\eqref{t123} are constant multiples of $t_1^{2d}(t_1+t_2)$ modulo
$(t_1+t_2)^2$. By \eqref{kkwp}, 
$$\Big \langle (d)\ \Big| -\sigma_1(F)\
\Big| \ (d) \Big \rangle_n = \frac{\gamma_{d,n}}{t_1^2} (t_1+t_2)\ \mod
(t_1+t_2)^2.$$
Hence,
we need only  
verify the equality (\ref{t123}) modulo $(t_1+t_2)^2$.

\subsection{Induction step: I}
Let $d>1$ and let $n> d$.
For the induction step, we will compute the invariant
\begin{equation}
  \label{inv_comp}
  \Big \langle  \left[\cI_{(d)}\right] \ \Big |
 -\sigma_1(F) \ \Big| \  \left[\cI_{(d-1,1)}\right]
\Big \rangle_n \,.
\end{equation}
Here,
 $\cI_{\lambda}$ denotes the monomial ideal corresponding to the
partition $\lambda$, and $[\cI_\lambda]$ denotes the
$T$-equivariant class of the associated fixed point in 
$\Hilb(\com^2,{|\lambda|})$.

The $T$-fixed point $[\cI_{\lambda}]$
corresponds to the Jack polynomial $$\bJ^\lambda\in\cF\otimes \Q[t_1,t_2].$$ 
For $-t_2/t_1=1$, the Jack polynomials specialize to the Schur functions.
Hence,
$$
\bJ^\lambda \equiv \frac{(-1)^{|\lambda|} \, |\lambda|!}{\dim \lambda}
\sum_\mu \chi^\lambda_\mu  \, t_1^{|\lambda|+\ell(\mu)}\, 
\bv \mu \brang  \mod t_1+t_2 \,,
$$
where $\dim \lambda$ is the dimension of the
representation $\lambda$ of the symmetric group and 
$\chi^\lambda_\mu$ is the evaluation of the
corresponding character on the conjugacy class $\mu$, see \cite{Mac}. 
In particular, the coefficient of $\bv d \brang$
in the expansion of both $\bJ^{(d)}$ and $\bJ^{(d-1,1)}$
is nonzero.

\begin{Lemma} We have \label{indone}
\begin{multline}
\label{JJ}
\lang \bJ^{(d)} \Bv \MM_\sigma-\MM_\sigma(0) \Bv  \bJ^{(d-1,1)} 
\rang_{\mathcal F}
\equiv \\
(-1)^d (t_1+t_2) \,\frac{t_1^{2d} \, (d!)^2}{d-1} 
\left(\frac{q}{1+q}+d\frac{(-q)^{d}}{1-(-q)^{d}}\right)
 \mod  (t_1+t_2)^2\,. 
\end{multline}
\end{Lemma}
Here, the operator $\MM_\sigma-\MM_\sigma(0)$ is formed by the 
terms of positive degree in $q$ in the operator $\MM_\sigma$.

The proof of Lemma \ref{indone} follows exactly the 
derivation of equation  (31) of \cite{qchs}. 
The differences between our operator $\MM_\sigma$ and 
the operator $\MM_D$ of \cite{qchs}  yield
constant functions orthogonal to the non-trivial
character $\chi^{(d-1,1)}$.

\subsection{Localization}
\subsubsection{Overview}
Our goal now is to reproduce the answer \eqref{JJ} by calculating the 
Donaldson-Thomas invariant
\begin{equation*}
  \lang \left[\cI_{(d)}\right]\Big| 
 -\sigma_1(F)\ \Big|  \left[\cI_{(d-1,1)}\right]
\rang_n \, .
\end{equation*}
By the rubber calculus \eqref{dx},
\begin{equation} \label{vvbb}
\lang \left[\cI_{(d)}\right]     \Big| 
 -\sigma_1(F)\ \Big|    \left[\cI_{(d-1,1)}\right]
\rang_n  = (n-d) \lang \left[\cI_{(d)}\right]  \Big|  
\left[\cI_{(d-1,1)}\right] \rang^\sim_n. 
\end{equation}
We calculate the right side of \eqref{vvbb} by
 $T$-equivariant localization on the
rubber moduli space 
$$R_{n,d}= I_n(R/R_0\cup R_\infty)^\sim,$$
see Section \ref{clrt}.

Since the $T$-fixed locus of the moduli space 
$R_{n,d}$ is proper, the
 virtual localization formula of \cite{gp} may be applied. 
However, since $R_{n,d}$
contains positive dimensional families of 
$T$-invariant ideal sheaves, a straightforward
application is difficult. 
Our strategy for  computing \eqref{vvbb}
uses comparisons to integrals in the quantum cohomology of the
Hilbert scheme of points $\text{Hilb}(\com^2,d)$.

\subsubsection{Skewers and twistors}
Consider the rubber moduli space $R_{n,d}$ for $n>d$.
Let 
$$[I] \in R^T_{n,d}$$
be a $T$-fixed ideal sheaf.
The ideal sheaf $I$ is defined on a 
rubber target fibered by $\com^2$ over a chain $C$ of
rational curves. The diagram below 
gives an example of a  subscheme associated
to a $T$-fixed ideal sheaf on a reducible rubber
target.

 \begin{center}
    \scalebox{0.7}{\includegraphics{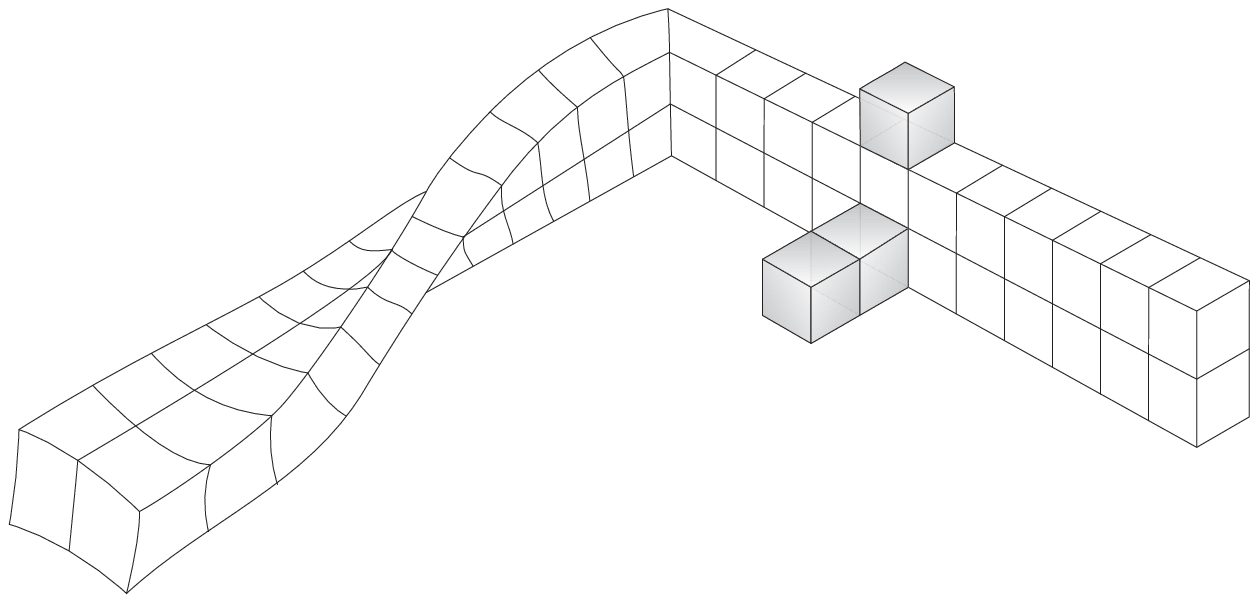}}
  \end{center}

Let $P\subset C$ be a rational component. The restriction
 $I_P$ of $I$
to the component
$\com^2\times P$ of the target rubber
is either a {\em skewer} or a {\em twistor}:
\begin{enumerate}
\item[(i)] A skewer is determined by an element  
$$[I_{sk}]\in I_*(\com^2\times \Pp,d)^T$$
where the factor $\Pp$ is rigid.
Certainly the skewer $I_{sk}$ is $T$-fixed on the rubber $\com^2 \times P$.
The component on the right of the above diagram is a skewer.

\item[(ii)]A twistor is obtained from a $T$-fixed element 
$$[f_{tw}]\in {M}_{0,\{0,\infty\}}\Big(\text{Hilb}\big((\com^2,d),*\big)
\Big)^T.$$
By pulling-back the universal ideal sheaf over the Hilbert scheme,
$f_{tw}$ determines a $T$-fixed ideal sheaf $I_{tw}$ on
the rubber $\com^2\times P$. The component on the left
of the diagram is a twistor.

The degree of the map $f_{tw}$ to the Hilbert scheme and the
Euler number $n_{tw}$ of $I_{tw}$ are related by
$$\text{deg}(f_{tw})= n_{tw}-d.$$

\end{enumerate}
Elementary considerations show
these disjoint constructions exhaust all $T$-fixed ideal sheaves
on the rubber $\com^2\times P$.

\subsubsection{Comparison}
We will calculate \eqref{vvbb} via a comparison result.
\begin{Proposition} \label{bcc} We have
$$
\lang \left[\cI_{(d)}\right]  \Big|  
\left[\cI_{(d-1,1)}\right] \rang^{ \sim}_n=
\lang \left[\cI_{(d)}\right]  \Big|  
\left[\cI_{(d-1,1)}\right] \rang^{\text{Hilb}(\com^2,d)}_{n-d}
\ \mod (t_1+t_2)^2.$$
\end{Proposition}
\noindent Proposition \ref{bcc} is proven in Section \ref{cvbn}.
As a first step, a simpler 
comparison result is obtained here.

Consider the open set of stable maps, 
$$U_{n,d} \subset \overline{M}_{0,\{0,\infty\}}\Big(
\text{Hilb}\big(\com^2,d),n-d\Big)$$
for which the domain is a {\em chain}
of rational curves.
The open set $U_{n,d}$ carries a $T$-equivariant Gromov-Witten
obstruction theory via restriction.
By pulling-back the universal ideal sheaf over the Hilbert
scheme (as in the twistor construction), we obtain an
open immersion
$$\iota: U_{n,d} \subset R_{n,d}.$$
Hence, $U_{n,d}$ also carries a $T$-equivariant Donaldson-Thomas
theory.

\begin{Lemma} The $T$-equivariant Gromov-Witten and Donaldson-Thomas 
obstruction 
theories of $U_{n,d}$ are canonically isomorphic.\label{ccfftt}
\end{Lemma}

\begin{proof} 
Let $$[\ f: C \rarr \text{Hilb}(\com^2,d) \ ] \in U_{n,d}$$
be a stable map, and let
  $$\iota([f]) = [I] \in R_{n,d}$$
be the associated ideal sheaf.
Since both obstruction theories
can be defined relative to the degenerations
of the domain,
we need only check the Gromov-Witten
complex associated to the cohomology 
$H^*(C,f^* T_{\text{Hilb}})$ 
matches the shift of the complex associated to $\text{Ext}_0(I,I)$. 
For a point $p\in C$, the tangent to the Hilbert scheme
may be viewed as $\text{Ext}_0^1(I_p,I_p)$. Moreover,
$$\text{Ext}_0^0(I_p,I_p)=\text{Ext}_0^2(I_p,I_p)=0.$$
The required matching
is then a straightforward application of the Leray spectral sequence.
\end{proof}

\subsubsection{Tangent representations}
Let $[I]\in R_{n,d}^T$ be a $T$-fixed ideal sheaf on a rubber target
fibered by $\com^2$ over a chain $C$ of rational curves.

The {\em special} points of $C$ consist of $0,\infty\in C$
and all the nodes.
Over the special points of $s \in 
C$, the ideal $I$  must correspond to $T$-fixed
points 
$[I_s]\in \text{Hilb}(\com^2,d)$.

Let $P \subset C$ be a component containing the special points
 $s,s'$ of $C$.
Fractional $T$-weight $w_{P,s}$ and $w_{P,s'}$ are defined by the 
$T$-representation of the
tangent space to $P$ at $s$ and $s'$:
\begin{enumerate}
\item[(i)] If $P\subset C$ corresponds to a skewer, then
$$[I_s]=[I_{s'}]$$
and
$w_{P,s}=w_{P,s'}=0.$
\item[(ii)] If
$P\subset C$ corresponds to a twistor, then
$$[I_s]\neq [I_{s'}]$$
and 
$w_{P,s},w_{P,s'}\neq 0$ modulo $(t_1+t_2)$, see \cite{qchs}.
\end{enumerate}
In fact, in the twistor case, if
$$[I_s]=[\cI_\mu], \ \ [I_{s'}]= [\cI_\nu],$$
then
$$w_{P,s}= w^{n_P-d}_{\mu,\nu}= - w^{n_P-d}_{\nu,\mu} =-w_{P,s'},$$
where $n_P$ is the Euler characteristic associated to $I_P$ and
$w^n_{\mu,\nu}$ is the universal function defined in Lemma 4 of
\cite{qchs} viewed here as a $T$-weight.

\subsubsection{Proof of Proposition  \ref{bcc}}\label{cvbn}
We index the $T$-equivariant localization contributions to the 
invariants 
$$\lang  [\cI_\mu] \ \Big| \ [\cI_\nu] \rang^\sim_{n}$$
by graphs following \cite{qchs}.
An {\em oriented chain} of Euler number $n$ is a graph 
$\Gamma_{\mu,\nu}=(V,v_1,v_2,\rho,E,S, \delta)$, 
\begin{enumerate}
\item[(i)] $V$ is a finite vertex set with distinguished elements $v_1$ and $v_2$,
\item[(ii)] $\rho: V \rightarrow {\mathcal P}(d)$, 
\item[(iii)] $E$ is a finite edge set,
\item[(iv)] $S\subset V$ ,
\item[(v)] $\delta:E\cup S \rightarrow {\mathbb Z}_{>0}$ is an assignment,
\end{enumerate}
satisfying the following conditions
\begin{enumerate}
\item[(a)] $\Gamma$ is a connected chain with initial vertex $v_1$ and final
vertex $v_2$,
\item[(b)] $\rho(v_1)=\cI_\mu,$ $\rho(v_2)=\cI_\nu$,
\item[(c)] if $v',v'' \in V$
are connected by an edge, then $\rho(v')\neq \rho(v'')$,
\item[(d)] if $v\in V\setminus S$ has edge valence 2 with neighbors $v',v''$, then
$$w^{\delta(e(v,v'))}_{\rho(v),\rho(v')} + w^{\delta(e(v,v''))}_{\rho(v),\rho(v'')}\neq 0
\ \mod t_1+t_2,$$
\item[(e)] $\sum_{e\in E} \delta(e) +\sum_{s\in S} \delta(s) =n+d(|E|+|S|-1).$
\end{enumerate}

Let $[I]\in R_{n,d}^T$ be $T$-fixed ideal sheaf on a rubber target
fibered by $\com^2$ over a chain $C$ of rational curves.
We associate an oriented chain, 
$$\Gamma_I=(V,v_1,v_2,\rho,E,S,\delta),$$ of Euler number
 $n$ to $I$ by the following
construction.
The vertex set 
$$V=V_1\cup V_2\cup V_3,$$ is a union of three disjoint subsets:
\begin{enumerate}
\item[(1)] $V_1$ is the set of maximal connected subcurves
 of skewer components of the base $C$,
\item[(2)] $V_2$ is the set of nodes of $s\in C$ for which the incident
components $P,P'\subset C$ are twistors and satisfy the {\em breaking
condition}
$$w_{P,s}+w_{P',s} \neq 0 \  \mod t_1+t_2,$$
\item[(3)] $V_3$ is subset of the marking $0,\infty \in C$
which lie on twistors.
\end{enumerate}
The markings $0,\infty \in C$ are 
associated to elements of the union $V_1\cup V_3$ ---
the markings  determine $v_1$ and $v_2$. 
The function $\rho$ is obtained from $I$.
Chains of unbroken twistors of $C$ link the vertices of $V$.
The edge set $E$  is determined by such chains. 
The set $S$ equals $V_1$.
The degree assignment $\delta$ is obtained from the 
Euler number of the restriction of $I$ to the associated components.

The oriented chain $\Gamma_I$ is easily seen to satisfy 
conditions (a)-(e). Condition (c) is a consequence of Lemma 5 of
\cite{qchs}.
The chain $\Gamma_I$ is invariant as $[I]$ 
varies in a connected component of the
$T$-fixed locus of $R_{n,d}^T$.

Let $G_n(\mu,\nu)$ denote the finite set of oriented chains trees 
$\Gamma_{\mu,\nu}$ of Euler
number $n$.
Let $$R_\Gamma \subset R^T_{n,d}$$
denote 
the substack of $T$-fixed maps corresponding to the tree 
$\Gamma\in G_n(\mu,\nu)$.
Let
$$\lang   [\cI_\mu]\ \Big| \ [\cI_\nu]
  \rang_{d}^{\Gamma\sim}$$
denote the localization contribution of $R_\Gamma$.
By the virtual localization formula \cite{gp},
\begin{equation}\label{kk12}
\lang [\cI_\mu] \ \Big| \ 
[\cI_\nu] \rang_d^{\Gamma\sim} = \int_{[R_\Gamma]^{vir}} 
\frac{\epsilon_0^*([\cI_\mu]) \cup
\epsilon_\infty^*([\cI_\nu])}{e({N}^{vir})}.
\end{equation}

By decomposing the Donaldson-Thomas obstruction theory, 
we can express the integral \eqref{kk12} in terms of descendent
skewer integrals corresponding to the vertices $S$ and
descendent twistor integrals corresponding to the edges $E$.

The edge integrals are exactly equal to associated
integrals in the Gromov-Witten theory of the
Hilbert scheme $\text{Hilb}(\com^2,d)$ by Lemma \ref{ccfftt}.
Moreover, the Hilbert scheme integrals which arise
for each edge compute a Gromov-Witten residue {\em {for
the $\TT$-action}}. Hence, the edge term
contributions to \eqref{kk12} are divisible by $(t_1+t_2)$.

We prove the skewer integrals are also divisible by $(t_1+t_2)$
by localizing a rubber integral:
\begin{multline*}
\lang [\cI_\gamma] \ \Big| \ \psi_0^a \psi_\infty^b \ \Big| 
\ [\cI_\gamma] \rang^\sim_n =
\lang [\cI_\gamma] \ \Big| \ \psi_0^a \psi_\infty^b \ \Big| 
\ [\cI_\gamma] \rang^{\Gamma^0_{\gamma,\gamma}\sim}\\ +
\sum_{\Gamma \in G_n(\gamma,\gamma)\setminus \{\Gamma^0\}} 
\lang [\cI_\gamma] \ \Big| \ \psi_0^a \psi_\infty^b \ \Big| 
\ [\cI_\gamma] \rang^{\Gamma \sim}.
\end{multline*}
Here, $\Gamma^0_{\gamma,\gamma}$ is the unique chain with a single skewer
vertex.
By Lemma \ref{rubvan}, the rubber integral on the left is divisible by
$(t_1+t_2)$. 
The $\Gamma^0$ contribution is the skewer integral of interest.
The second summand on the right is expressed in terms of 
skewer integrals of {\em lower} Euler number and edge integrals.
The former are divisible by $(t_1+t_2)$ by induction. 
We conclude the skewer integrals are divisible by $(t_1+t_2)$.

We now specialize to the localization analysis of the rubber integral,
$$\lang  [\cI_{(d)}] \ \Big| \ [\cI_{(d-1,d)}] \rang^\sim_{n} =
\sum_{\Gamma \in G_n(\mu,\nu)} 
\lang [\cI_{(d)}] \ \Big|  
\ [\cI_{(d-1,1)}] \rang^{\Gamma \sim}.$$
Assume $\mu\neq \nu$.
Since both skewer vertices and edges contribute factors
of $(t_1+t_2)$,
$$\lang  [\cI_{(d)}] \ \Big| \ [\cI_{(d-1,1)}] \rang^\sim_{n} = 
\lang [\cI_{(d)}] \ \Big|  
\ [\cI_{(d-1,1)}] \rang^{\Gamma^0_{\mu,\nu}\sim} \ \mod t_1+t_2,$$
where $\Gamma^0_{\mu,\nu}$ is the unique single edged chain
with $S=\emptyset$. The Proposition \ref{bcc} is then a
consequence of Lemma \ref{ccfftt}.\qed

\subsection{Induction step: II}
By Lemma \ref{indone}, Proposition \ref{bcc}, and the Hilbert scheme
calculation of \cite{qchs}, we obtain
$$
\lang \bJ^{(d)} \Bv \MM_\sigma-\MM_\sigma(0) \Bv  \bJ^{(d-1,1)} 
\rang_{\mathcal F, n-d} = \lang \left[\cI_{(d)}\right] \  \Big| -\sigma_1(F)\
\Big| \  
\left[\cI_{(d-1,1)}\right] \rang^{ \sim}_n,$$
completing the proof of Proposition \ref{T12}.\qed


\section{The level $(0,0)$ theory}
\subsection{The operator $\MM_D$}
Let $D$ denote the insertion of the relative condition
$- (2, 1^{d-2})$
in the local Donaldson-Thomas theory of curves.
Let the operator ${\MM}_D$ on
Fock space be defined by matrix elements
\begin{eqnarray*}
\Big \langle \mu \ \Big| \ {\MM}_D \ \Big| \nu \Big \rangle_{\mathcal F} & = &
\mathsf{DT}(0|0,0)_{\mu,D,\nu}\\
& = & 
-\mathsf{DT}(0|0,0)_{\mu,(2, 1^{d-2}),\nu}.
\end{eqnarray*}
By definition, the insertion $D$ vanishes in degrees $d=0,1$.

\begin{Proposition}\label{T22} 
$
    {\MM}_D = \MM - \frac{t_1+t_2}{2} \, \frac{(-q)+1}{(-q)-1} \, 
|\,\cdot\,| \,\,.
$
\end{Proposition}

\begin{proof} By applying the degeneration formula to the
definition of $\MM_\sigma$, we obtain
\begin{eqnarray*}
\Big \langle \mu \ \Big| {\MM}_\sigma  \Big| \nu \Big \rangle_{\mathcal F} 
& = &
q^{-d} \Big \langle \mu \ \Big| -\sigma_1(F)\ \Big|\ \nu \Big \rangle\\
& = & \sum_\gamma q^{-d} \mathsf{Z}(0|0,0)_{\mu,\gamma,\nu} 
\bigtriangleup_d(\gamma,\gamma) 
q^{-d}\Big\langle \gamma \ \Big | -\sigma_1(F) \Big\rangle \\
& = & \sum_\gamma \mathsf {DT}(0|0,0)_{\mu,\gamma,\nu}
\bigtriangleup_d(\gamma,\gamma) q^{-d}
\frac{\Big\langle \gamma \ \Big | -\sigma_1(F) \Big\rangle}{
\Big\langle \emptyset \ \Big |  \Big\rangle}
\end{eqnarray*}
By equation \ref{van7},
$$\Big\langle \gamma \ \Big | -\sigma_1(F) \Big\rangle=0$$
if $\ell(\gamma)<d-1$. Hence,
there are only two nonvanishing terms in the sum over the partition
$\gamma$:
\begin{eqnarray*}
\Big \langle \mu \ \Big|  {\MM}_\sigma  \Big| \nu \Big \rangle_{\mathcal F}
 & = &\  \
\mathsf {DT}(0|0,0)_{\mu,(1^d),\nu} \frac{d!(t_1t_2)^d}{q^d}
\frac{\langle (1^d)  | -\sigma_1(F) \rangle}{
\langle \emptyset \  |  \rangle}\\
& & -
\mathsf {DT}(0|0,0)_{\mu,(2,1^{d-2}),\nu} \frac{2(d-2)!(t_1t_2)^{d-1}}{q^d}
\frac{\langle (2,1^{d-2})  | -\sigma_1(F) \rangle}{
\langle \emptyset \  |  \rangle}
\end{eqnarray*}

By Lemma \ref{lcap} for the cap $\mathsf{DT}(0|0,0)_\lambda$,
the insertion $(1^d)$ can be freely added or removed in the
local Donaldson-Thomas theory of curves.
Hence,
\begin{eqnarray*}
\mathsf {DT}(0|0,0)_{\mu,(1^d),\nu}& =& \mathsf {DT}(0|0,0)_{\mu,\nu}\\
                                    & =& 
\Big\langle \mu\ \Big| \ \nu \Big \rangle_{\mathcal F},
\end{eqnarray*}
Similarly,
 \begin{eqnarray*}
q^{-d}\frac{\Big\langle (1^d) \  \Big| -\sigma_1(F) \Big \rangle}{
\Big \langle \emptyset \  \Big|  \Big\rangle} & = & 
q^{-d}\frac{\Big\langle (1^d) \  \Big| -\sigma_1(F)\ \Big|\ (1^d) \Big\rangle}{
\Big\langle \emptyset \  \Big|\ \emptyset \Big\rangle}\\
& = & \Big \langle (1^d) \ \Big| \ {\MM}_\sigma \ \Big| (1^d) \Big \rangle\\
& = & \frac{(t_1+t_2)}{d!(t_1t_2)^d} \left( 
\frac{d}{2}\frac{(-q)+1}{(-q)-1}
-\Phi(q)\right),
\end{eqnarray*}
where the last equality is obtain from Proposition \ref{T12}.
Finally,
\begin{eqnarray*}
q^{-d}\frac{\Big\langle (2,1^{d-2}) \  \Big| -\sigma_1(F) \Big \rangle}{
\Big \langle \emptyset \  \Big|  \Big\rangle} & = & 
q^{-d}\frac{\Big\langle (2,1^{d-2}) 
\  \Big| -\sigma_1(F)\ \Big|\ (1^d) \Big\rangle}{
\Big\langle \emptyset \  \Big|\ \emptyset \Big\rangle}\\
& = & \Big \langle (2,1^{d-2}) 
\ \Big| \ {\MM}_\sigma \ \Big|\ (1^d) \Big \rangle_{\mathcal F}\\
& = & -\Big\langle (2,1^{d-2}) \ \Big| \ (2,1^{d-2}) \Big \rangle_{\mathcal F}.
\end{eqnarray*}
We conclude
\begin{eqnarray*}
\Big \langle \mu \ \Big| \ {\MM}_D \ \Big| \nu \Big \rangle_{\mathcal F} 
& = & -
\mathsf {DT}(0|0,0)_{\mu,(2,1^{d-2}),\nu} \\
& = & 
\Big \langle \mu \ \Big| \ {\MM}_\sigma \ \Big| \nu \Big \rangle_{\mathcal F}
-(t_1+t_2)
\Big\langle \mu\ \Big| \ \nu \Big \rangle_{\mathcal F}
\left( \frac{d}{2}\frac{(-q)+1}{(-q)-1}
- {\Phi(q)}\right) \\
& = & \Big \langle \mu \ \Big| \ 
 \MM - \frac{t_1+t_2}{2} \, \frac{(-q)+1}{(-q)-1}\, 
|\,\cdot\,|
 \ \Big| \nu \Big \rangle_{\mathcal F}\, ,
\end{eqnarray*}
where the last equality is obtained from Proposition \ref{T12}.
\end{proof}

\subsection{Proof of Theorem \ref{ccc} in level $(0,0)$}

\subsubsection{Starred series}
The terminology of \cite{bp} is convenient for discussing
the GW/DT correspondence. On the Gromov-Witten side, let
$${\mathsf{GW}}^*
(g|  k_{1},k_{2})_{\lambda ^{1}\dots \lambda ^{r}} =
(-iu)^{d (2-2g+k_{1}+k_{2})-\delta } \ 
{\mathsf Z}'(N)_{ \lambda^1 \dots \lambda^r},$$
where $N$ is rank 2 bundle of level $(k_1,k_2)$ on a genus $g$
curve and
$$\delta = \sum_{i=1}^r (d - \ell(\lambda^i)).$$
On the Donaldson-Thomas side, let
\begin{eqnarray*}
\mathsf{DT}^*(g| k_1,k_2)_{\lambda^1\dots\lambda^r} & = & 
(-q)^{-\frac{d}{2}(2-2g+k_1+k_2)} \ {\mathsf Z}'_{DT}
(N)_{\lambda^1 \dots \lambda^r} \\ & =& 
(-1)^{-d(1-g)}(-q)^{-\frac{d}{2}(k_1+k_2)} {\mathsf{DT}}(g|k_1,k_2).
\end{eqnarray*}
Theorem \ref{ccc} of the  Gromov-Witten/Donaldson-Thomas 
correspondence for local curves can be
restated as the equality
\begin{equation}\label{kk231}
\mathsf{GW}^*(g|  k_1,k_2)_{\lambda^1\dots\lambda^r} =
\mathsf{DT}^*(g|  k_1,k_2)_{\lambda^1\dots\lambda^r},
\end{equation}
after the variable change $e^{iu}=-q$,
\subsubsection{TQFT}
Slightly altered metrics are defined for raising the indices of
the starred series:
\begin{eqnarray*}
{\mathsf{GW}}^* (g| k_{1},k_{2})_{\mu ^{1}\dots \mu ^{s}}^{\nu
^{1}\dots \nu ^{t}} & =& \left(\prod _{i=1}^{t} \zz(\nu ^{i})
(-t_{1}t_{2})^{l (\nu ^{i})} \right) {\mathsf{GW}}^* (g|
k_{1},k_{2})_{\mu
^{1}\dots \mu ^{s}\nu ^{1}\dots \nu ^{t}},
\\
{\mathsf{DT}}^* (g| k_{1},k_{2})_{\mu ^{1}\dots \mu ^{s}}^{\nu
^{1}\dots \nu ^{t}}& =&  \left(\prod _{i=1}^{t} \zz(\nu ^{i})
(-t_{1}t_{2})^{l (\nu ^{i})} \right) {\mathsf{DT}}^* (g|
k_{1},k_{2})_{\mu
^{1}\dots \mu ^{s}\nu ^{1}\dots \nu ^{t}}.
\end{eqnarray*}
With respect to the above metrics, the starred partition
function satisfy the same degeneration rules as their
unstarred counterparts:
\begin{eqnarray*}
{\mathsf {GW}}^*(g|k_1,k_2)_{\mu^1,\dots,\mu^s}^{\nu^1,\dots,\nu^t} &
 = &
\sum_\gamma
{\mathsf{GW}}^*(g'|k'_1,k'_2)_{\mu^1,\dots,\mu^{s}}^\gamma  
{\mathsf{GW}}^*(g''|k''_1,k''_2)^{\nu^{1},\dots,\nu^t}_{\gamma}\, ,\\
{\mathsf {DT}}^*(g|k_1,k_2)_{\mu^1,\dots,\mu^s}^{\nu^1,\dots,\nu^t}&  =&
\sum_\gamma
{\mathsf{DT}}^*(g'|k'_1,k'_2)_{\mu^1,\dots,\mu^{s}}^\gamma  
{\mathsf{DT}}^*(g''|k''_1,k''_2)^{\nu^{1},\dots,\nu^t}_{\gamma}\, ,
\end{eqnarray*}
where $g=g'+g''$, and $k_i=k_i'+k_i''$, and
\begin{eqnarray*}
{\mathsf{GW}}^*(g|k_1,k_2)_{\mu^1,\dots,\mu^s} & = &
\sum_\gamma {\mathsf{GW}}^*(g-1|k_1,k_2)^\gamma_{\mu^1,\dots,\mu^s,
\gamma}\, , \\
{\mathsf{DT}}^*(g|k_1,k_2)_{\mu^1,\dots,\mu^s} & =&
\sum_\gamma {\mathsf{DT}}^*(g-1|k_1,k_2)^\gamma_{\mu^1,\dots,\mu^s,
\gamma}\, .
\end{eqnarray*}
Hence, tensor functors
\[
{\mathbf{GW^*}},{\mathbf{DT^*}}
:2\mathbf{Cob}^{L_{1},L_{2}}\to R\mathbf{mod}.
\]
can be defined just as before.

\subsubsection{Matching in level $(0,0)$}
Using the TQFT structure, to prove Theorem \ref{ccc}
in level $(0,0)$, we must establish the following three equalities:
\begin{eqnarray*}
{\mathsf{GW}}^*(0|0,0)_\lambda & = & 
{\mathsf{DT}}^*(0|0,0)_\lambda, \\
{\mathsf{GW}}^*(0|0,0)_{\lambda \mu} & = & 
{\mathsf{DT}}^*(0|0,0)_{\lambda \mu}, \\
{\mathsf{GW}}^*(0|0,0)_{\lambda \mu \nu} & = & 
{\mathsf{DT}}^*(0|0,0)_{\lambda\mu \nu},
\end{eqnarray*}
corresponding respectively to the cap, the tube, and the pair
of pants.

The matching of the level $(0,0)$ cap is a consequence of
Lemma 6.2 of \cite{bp} for the Gromov-Witten side and
Lemma \ref{lcap} of Section \ref{ddd111} for the Donaldson-Thomas side.
Similarly, the level $(0,0)$ tube matching is
a consequence of Lemma 6.1 of \cite{bp} and Lemma \ref{feww} of
Section \ref{tqft}.

The pair of pants matching in level $(0,0)$ is more subtle.
The main result of the Appendix of \cite{bp} is the
unique determination of the level $(0,0)$ TQFT for the local 
Gromov-Witten theory of curves by the cap, the tube,
and the set of series
$${\mathsf{GW}^*}(0|0,0)_{\mu,(2,1^{d-2}),\nu}$$ for all $\mu$ and $\nu$.
Since the cap and the tube have been shown to match,
the equality 
\begin{equation}\label{xxgg6}
{\mathsf{GW}^*}(0|0,0)_{\mu,(2,1^{d-2}),\nu} =
{\mathsf{DT}^*}(0|0,0)_{\mu,(2,1^{d-2}),\nu}
\end{equation}
suffices to complete the matching in level $(0,0)$.
Equality \eqref{xxgg6} is a consequence of A.3 of \cite{bp} and
Proposition \ref{T22} above. \qed

\section{The cap of level $(-1,0)$}
\label{xxx}

\subsection{TQFT}
By Theorem 4.1 of \cite{bp}, the proof of Theorem \ref{ccc} for
all levels now
requires only the equality
\begin{equation}\label{kkk231}
\mathsf{GW}^*(0|-1,0)_{\lambda} =
\mathsf{DT}^*(0|-1,0)_{\lambda}.
\end{equation}
The Gromov-Witten side was calculated in \cite{brp,bp}.
By Lemma 6.3 and Section 6.4.1 of \cite{bp},
$$\mathsf{GW}^*(0|-1,0)_{\lambda} = (-q)^{-d/2} t_2^{-\ell(\lambda)}
\frac{(-1)^{d-\ell(\lambda)}}
{\zz(\lambda)} t_2^{-\ell(\lambda)}
\prod_{i=1}^{\ell(\lambda)} \frac{1}{1-(-q)^{-\lambda_i}}$$
In order to prove \eqref{kkk231}, we must find the following
evaluation:
\begin{eqnarray*}
\mathsf{DT}(0|-1,0)_\lambda & =&
 (-1)^d (-q)^{-d/2} \mathsf{DT}^*(0|-1,0)_\lambda \\
& = &
 {t_2^{-\ell(\lambda)}}
\frac{(-1)^{d-\ell(\lambda)}}
{\zz(\lambda)} q^{-d}
\prod_{i=1}^{\ell(\lambda)} \frac{1}{1-(-q)^{-\lambda_i}} \\
& = & \frac{t_2^{-\ell(\lambda)}}
{\zz(\lambda)} 
\prod_{i=1}^{\ell(\lambda)} \frac{1}{1-(-q)^{\lambda_i}}.
\end{eqnarray*}

\subsection{${\mathbf{T}}$-action}
Let $T$ be the standard 2-dimensional torus action on
the bundle 
$$
N = \cO(-1) \oplus \cO 
$$
over $\Pp$ with scaling weights $t_1$ and $t_2$ on the
factors $\cO(-1)$ and $\cO$ respectively.
Let the 1-dimensional torus $S$ act with weights
$s,-s$ at the fixed points $$0,\infty \in \Pp$$
and weights $(-s,0)$ on the fiber of $N_\infty$.
We will consider the full
$${\mathbf T}= S \times T $$
action on $N$.

\subsection{Independence}
Let $L\subset N$ denote the ${\mathbf T}$-equivariant
divisor determined by the summand
$\cO(-1)$, and let
$$[L] \in A^1_{\mathbf T}(N, {\mathbb Q})$$
be the associated class.
 Let $$\lambda([L])= \{ \lambda_1([L]), \ldots,
\lambda_{\ell(\lambda)}([L])\} $$ be a weighted
partition of $d$.
The reduced ${\mathbf T}$-equivariant Donaldson-Thomas
residue integral relative to $N_\infty$,
\begin{equation*}
\Big \langle\   \Big| \ \lambda([L]) \Big \rangle'_{(-1,0)}=
\frac{\Big \langle \  \Big| \ \lambda([L]) \Big \rangle_d^{(-1,0)}}
{\Big \langle\   \Big| \emptyset \Big \rangle^{(-1,0)}_0},
\end{equation*}
is of degree 0 in the equivariant parameters $s$, $t_1$, and $t_2$.

\begin{Lemma}\label{indd45}
$\Big \langle\   \Big| \ \lambda([L]) \Big \rangle'_{(-1,0)}$
 is 
{\em independent} of $s$, $t_1$, and $t_2$.
\end{Lemma}

\begin{proof}
Let $N \subset \overline{N}$
be a ${\mathbf T}$-equivariant compactification, and
let 
$$[P]\in H_2(\overline{N},{\mathbb Z})$$
be the push-forward to $\overline{N}$ 
of the class of the zero section $P\subset N$.
Consider the restricted moduli space of ideal sheaves
$$R=\epsilon^{-1}_{\infty} (C_{\lambda[L]}) \subset
I_n(\overline{N}/\overline{N}_\infty,d[P]),$$
following the notation of Section \ref{relgeom}.
Since the line bundle $\cO(-1)$ over $\Pp$ has no nontrivial multisections,
the elements of 
$$[I_Z] \in R_n$$ for which
the
entire 1-dimension support of $Z$ lies
on $P$ determine an open and closed 
${\mathbf T}$-equivariant substack $R^0_n\subset R_n$

Dimension 0 integrals over $R_n^0$ are certainly independent of
the equivariant parameters $s$, $t_1$, and $t_2$.
By localization,
$$
\Big \langle\   \Big| \ \lambda([L]) \Big \rangle'_{(-1,0)}=
\frac{\sum_n q^n \int_{[R_n^0]^{vir}} 1}{\sum_n q^n
\int_{[I_n(\overline{N}/\overline{N}_\infty,0)]^{vir}} 1}.$$
Since the denominator on the right is also independent of
the equivariant parameters, the Lemma is proven. 
\end{proof}

\subsection{Localization}

The ${\mathbf T}$-equivariant
virtual localization formula for the
series 
\begin{equation}\label{vv23}
\Big\langle\   \Big| \ \lambda([L])  \Big\rangle'_{(-1,0)}
\end{equation}
involves an edge summation 
over all ${T}$-fixed points $\cI_\mu$ of the 
Hilbert scheme $\text{Hilb}(N_\infty,d)$. See
\cite{mnop1,mnop2} for a discussion of localization in
relative Donaldson-Thomas theory.

We  
orient the partition $\mu$ so that the rows of 
the associated Young diagram extend in the $\cO$ direction. 
Define $n(\mu)$ by a summation over rows:
\begin{equation}
  \label{nl}
  n(\mu) = \sum_{i=1}^{\ell(\mu)} (i-1)\mu_i \,. 
\end{equation}
With our orientation conventions, 
$d+n(\mu)$ is the Euler characteristic of 
a pure edge with profile $\mu$.  

By an application of the  virtual localization formula, we find 
the series  \eqref{vv23} equals
\begin{equation}\label{localiz}  
 \sum_{|\mu|=d} 
\frac{\bW(\mu,\emptyset,\emptyset)|_{s,t_1-s,t_2}}
{\bW(\emptyset,\emptyset,\emptyset)|_{s,t_1-s,t_2}}
  \cdot 
q^{n(\mu)}  
\bE^{(-1,0)}(\mu) 
\cdot 
\frac{\Big \langle \left[\cI_\mu \right]    \Big| \frac{1}{-s-\psi_\infty}
\ \Big| \lambda \Big \rangle^\sim}
{\Big \langle \emptyset    \Big| \frac{1}{-s-\psi_\infty}
\ \Big| \emptyset \Big \rangle^\sim}
\ t_2^{\ell(\lambda)}. 
\end{equation}
The terms   
$\bW$ and  $\bE^{(-1,0)}$
are respectively
the 
equivariant vertex and the
equivariant edge weight \cite{mnop1,mnop2}. The rubber integral
series in the numerator is
$$
\Big \langle  \left[\cI_\mu \right] \ \Big| \frac{1}{-s-\psi_\infty}
\ \Big| \lambda \Big \rangle^\sim =
q^d \Big\langle \left[\cI_\mu \right] \ \Big|  \lambda \Big \rangle_{d,d}
+ \sum_{n>d} q^{n} 
\Big \langle   \left[\cI_\mu \right]
\ \Big| \frac{1}{-s-\psi_\infty}
\ \Big| \lambda \Big \rangle^\sim_{n,d}.$$
The denominator series has a parallel definition.

By Lemma \ref{indd45}, the series \eqref{vv23}
is independent $s$, $t_1$, and $t_2$. Hence
the localization formula can be evaluated after specialization of 
the equivariant parameters.

\subsection{Rubber}
We evaluate the localization formula \eqref{vv23} after the 
specialization
\begin{equation} \label{sp234}
t_1+t_2=0.
\end{equation}
By the vanishing of Lemma \ref{rubvan}, the rubber integrals on
the right side simplify to
$$ 
\frac{\Big \langle \left[\cI_\mu \right]    \Big| \frac{1}{-s-\psi_\infty}
\ \Big| \lambda \Big \rangle^\sim}
{\Big \langle \emptyset    \Big| \frac{1}{-s-\psi_\infty}
\ \Big| \emptyset \Big \rangle^\sim} =
q^d \Big\langle \left[\cI_\mu \right] \ \Big| \ \lambda \Big \rangle_{d,d}.$$
The matrix element on the right is 
 the  equivariant intersection form of the classical
cohomology of $\text{Hilb}(N_\infty,d)$. 

For 
opposite weights \eqref{sp234}, 
the Hilbert scheme
intersections reduce to characters of the symmetric 
group, see for example \cite{vass}. We find 
\begin{equation}
  \label{inn_pr} 
\Big\langle \left[\cI_\mu \right] \ \Big| \ \lambda \Big \rangle_{d,d}
= 
\frac{t_2^{d-\ell(\lambda)}}{\zz(\lambda)} \, \chi^\mu_\lambda
 \, \prod_{\square\in \mu} h(\square) \,.
\end{equation}
Here, $h(\square)$ denotes the hook length. 
A similar calculation can be found in \cite{qchs}.

\subsection{Edge}

The edge term of the localization formula \eqref{localiz} 
also simplifies after the
specialization \eqref{vv23}.

We recall the formula for the edge term 
adapted to the level $(-1,0)$ geometry \cite{mnop1}.
Given a partition $\mu$, form
the following polynomials
\begin{eqnarray*}
Q_\mu(x_1,x_2)& =& \sum_{(i,j)\in \mu} x_1^{i} \, x_2^{j}\,, \\
F_\mu(x_1,x_2)& =&  - Q_\mu - \frac{\barQ_\mu}{x_1 x_2} +
Q_\mu \barQ_\mu \frac{(1-x_1)(1-x_2)}{x_1 x_2}.
\end{eqnarray*}
The sum in the first definition
 is over the interior corners of the squares of the 
 Young diagram of $\mu$ --- the corners closest to the origin. Also,
$$
\barQ_\mu(x_1,x_2) = Q_\mu(x_1^{-1},x_2^{-1}) \,.
$$
The rational function 
$$
E_\mu = \frac{F_\mu(x_1,x_2)}{x_0-1} + 
\frac{F_\mu(x_1x_0,x_2)}{x_0^{-1}-1}
$$
is readily seen to be a Laurent polynomial in the variables
$x_i$. The edge weight 
$\bE(\mu)$ 
is obtained from the following transformation: 
\begin{equation}
  \label{transf_E}
E_\mu = \sum_k a_k \, x_0^{k_0} x_1^{k_1} x_2^{k_2}   \mapsto 
\bE(\mu) = \prod_k (k_0 s+ k_1 (t_1-s) + k_2 t_2)^{-a_k} \,.
\end{equation}
Setting $t_1+t_2=0$ is equivalent to substituting 
$$
x_0 = (x_1 x_2)^{-1} 
$$
in the above formulas. 

\begin{Lemma} \label{arms&legs} We have 
  \begin{equation*}
    F_\mu(x_1,x_2) = - \sum_{\square\in \mu} 
\left(x_1^{l(\square)}x_2^{-a(\square)-1} + x_1^{-l(\square)-1}
x_2^{a(\square)}\right)\,,
  \end{equation*}
where $a(\square)$ and $l(\square)$ denote the arm-lengths and leg-length of 
a square in a diagram (number of squares to the right and below $\square$, 
respectively). 
\end{Lemma}

\begin{proof}
The polynomial $F_\mu$ is, up to sign, the character of the scaling
torus
action on the tangent space at $\left[\cI_\mu \right]$ to the
Hilbert scheme of points in ${\mathbb C}^2$. 
The exponents in the formula are well known to be the weights of tangent
action, see for example \cite{haiman}. 
\end{proof}

Define an auxiliary function 
$$
G(x_1,x_2) = - \sum_{\square\in \mu}x_1^{l(\square)}
 x_2^{-a(\square)-1} 
$$
for which $F = G + \barG/(x_1 x_2)$. We compute
\begin{equation}
  \label{E_}
E_\mu\big|_{x_3=(x_1 x_2)^{-1}} = - G|_{x_1=x_2^{-1}} + 
\frac{G-G|_{x_1=x_2^{-1}}}{(x_1 x_2)^{-1} - 1} 
 - \frac{\barG-\barG|_{x_2=x_1^{-1}}}{(x_1 x_2) - 1}  \,.
\end{equation}
All three terms of \eqref{E_} are Laurent polynomials.
The 
third term is minus bar of the second one. 
By the transformation \eqref{transf_E}, 
the factors corresponding to the second and third terms of \eqref{E_}  
cancel up to a sign. This sign is the parity of the overall number 
of monomials in the second term in \eqref{E_}  equal to
$$
\sum_{\square\in \mu} l(\mu) = n(\mu) \,.
$$
We have proven the following result:
  \begin{equation*}
   \bE(\mu)^{(-1,0)}|_{s,t_1-s,t_2} = (-1)^{n(\mu)} \, t_2^{-d} \, 
\prod_{\square\in \mu} h(\square)^{-1} \,. 
  \end{equation*}

\subsection{Vertex}\label{ver12}
The ${\mathbf T}$-equivariant tangent weights of $N$
at the fixed point over $0\in \Pp$ 
satisfy the Calabi-Yau condition after
the specialization \eqref{sp234}:
$$s + (t_1-s) + t_2 =0.$$
The vertex $\bW(\mu,\emptyset,\emptyset)$ has a rather
simple evaluation in the Calabi-Yau case \cite{mnop1,ORV}:
$$\bW(\mu,\emptyset,\emptyset)|_{s,t_1-s,t_2} =
\sum_\pi (-q)^{|\pi|}.$$
The sum is over all 3-dimensional partitions $\pi$ with
a single infinite leg in the $s$ direction asymptotic 
to $\mu$. Here, $|\pi|$ is the renormalized
volume --- the number of boxes remaining after the
infinite leg is removed. In particular,
$$\bW(\emptyset,\emptyset,\emptyset)|_{s,t_1-s,t_2} =
M(-q),$$
a specialization of \eqref{fwq}.

Evaluation of the Calabi-Yau vertex is reduced
to the enumeration of 3-dimensional partition. The enumeration
for the 1-legged vertex is solved in \cite{OR}, 
\begin{equation}\label{WF}
\frac{\bW(\mu,\emptyset,\emptyset)|_{s,t_1-s,t_2}}
{\bW(\emptyset,\emptyset,\emptyset)|_{s,t_1-s,t_2}} 
=
\prod_{\square \in \mu} \frac{1}{1-(-q)^{h(\square)}} \,, 
\end{equation}
where the product is over all squares in the Young diagram of $\mu$.

The origin of hook lengths
in \eqref{WF} is the following classical formula 
\cite{Mac} for 
the value of the Schur function $s_\mu$ at the 
point $(1,q,q^2,\dots)$, 
\begin{equation}
  \label{s(q)}
  s_\mu(1,q,q^2,\dots) = q^{n(\mu)} \, \prod_{\square \in \mu} 
\frac{1}{1-q^{h(\square)}} \,.
\end{equation}

\subsection{Evaluation}

Putting all pieces of the localization
formula \eqref{localiz} together, we find:
\begin{equation*}
  \Big \langle \lambda([L])  \Big| \Big \rangle'_{(0,-1)} = 
\frac{q^d}{\zz(\lambda)} \sum_\mu \chi^\mu_\lambda \, 
s_\mu(1,-q,(-q)^2,(-q)^3,\dots) \,.
\end{equation*}
By a classical formula in the theory of symmetric functions,
$\sum_\mu \chi^\mu_\lambda \, s_\mu$ equals the power sum
symmetric function $p_\lambda$, see \cite{Mac}. Since
$$
p_k(1,-q,(-q)^2,(-q)^3,\dots) = \frac{1}{1-(-q)^{k}}\,,
$$
we obtain the following result.

\begin{Proposition}\label{dddr}
  \begin{equation}
    \Big \langle \lambda([L])  \Big| \Big \rangle'_{(0,-1)} = \frac{q^d}{\zz(\mu)} 
\, \prod_{i=1}^{\ell(\lambda)} \frac{1}{1-(-q)^{\lambda_i}}
  \end{equation}
\end{Proposition}

\subsection{Proof of Theorem \ref{ccc}}
By definition, we find
$$\Big \langle \mu([L])  \Big| \Big \rangle'_{(0,-1)} = 
t_2^{\ell(\lambda)} q^d\  \mathsf{DT}(0|-1,0).$$
Hence, by Proposition \ref{dddr}, 
$$\mathsf{DT}(0|-1,0)=
\frac{t_2^{-\ell(\lambda)}}{\zz(\mu)} 
\, \prod_{i=1}^{\ell(\lambda)} \frac{1}{1-(-q)^{\lambda_i}}.$$
The matching \eqref{kkk231} of the $(-1,0)$ cap is established, and
the proof of Theorem \ref{ccc} is complete. \qed

\subsection{Proof of Theorem \ref{bbb}}
By Theorem \ref{ccc}, the rationality of
the Donaldson-Thomas series for local curves is a direct consequence
of Theorem 6.4 of \cite{bp}.
The proof of Theorem \ref{bbb} is complete. \qed

\section{The 1-legged vertex}
\subsection{Overview}
The localization formula for the level $(-1,0)$ cap 
together with a differential equation for rubber 
integrals provides an effective determination of the 1-legged
equivariant vertex.

\subsection{Differential equation}
Consider the following 
Donaldson-Thomas rubber descendent series:
$$
\Big \langle \mu \ \Big| \frac{1}{1-\psi_\infty}
\ \Big| \nu \Big \rangle^\sim =
q^d \Big\langle \mu \ \Big| \ \nu \Big \rangle_{d,d}
+ \sum_{n>d} q^{n} 
\Big \langle \mu \ \Big| \frac{1}{1-\psi_\infty}
\ \Big| \nu \Big \rangle^\sim_{n,d}.$$
Define a operator ${\mathsf S}$ on Fock space
by the matrix elements
\begin{equation}\label{lqde}
\Big\langle \mu \ \Big|
{\mathsf S} \ \Big| \nu \Big \rangle_{\mathcal F}
 = q^{-d} M(-q)^{t_1+t_2} \Big \langle \mu \ \Big| \frac{1}{1-\psi_\infty}
\ \Big| \nu \Big \rangle^\sim.
\end{equation}
By the rubber calculus relation \eqref{dx},
$$
q\frac{d}{dq} q^{-d}
\Big\langle \mu \ \Big| \frac{1}{1-\psi_\infty} \ \Big| \nu \Big \rangle
^\sim = 
 q^{-d}
\Big\langle \mu \ \Big| \frac{-\sigma_1(F)}
{1-\psi_\infty} \ \Big| \nu \Big \rangle
  -q^{-d}
\Big\langle \mu \ \Big| \frac{1}{1-\psi_\infty} \ \Big| D\cdot
\nu \Big \rangle^\sim.
$$
By the topological recursion relation of Section \ref{clrt},
$$
q^{-d}
\Big\langle \mu \ \Big| \frac{-\sigma_1(F)}{1-\psi_\infty} \ \Big| \nu \Big \rangle
= 
q^{-2d}
\sum_{\eta} 
\Big\langle \mu \ \Big| {-\sigma_1(F)} \ \Big| \eta \Big \rangle
\bigtriangleup_d(\eta,\eta) 
\Big\langle \eta \ \Big| \frac{1}{1-\psi_\infty} \ \Big| \nu \Big 
\rangle^\sim.$$
Together with Proposition \ref{T12}, we conclude
\begin{equation} \label{xxqq}
q\frac{d}{dq} {\mathsf S}
 = \MM {\mathsf S} - {\mathsf S} \MM(0),
\end{equation}
where $\MM(0)$ denotes the $q$-constant terms of $\MM$.
The series $\Phi(d)$ drops out of right side of \eqref{xxqq}.

The differential equation \eqref{xxqq} for
rubber descendents is almost identical to the quantum differential
equation for the Hilbert scheme of points of the plane \cite{qchs}.

\begin{Lemma} \label{cllw}
$\Big \langle \emptyset \ \Big| \frac{1}{1-\psi_\infty}
\ \Big| \emptyset \Big \rangle^\sim =
M(-q)^{-(t_1+t_2)}.$
\end{Lemma}

\begin{proof} The differential equation \eqref{xxqq} takes a
simple form in degree 0:
$$q \frac{d}{dq} \Big \langle \emptyset \ \Big| {\mathsf S}
\ \Big| \emptyset \Big \rangle_{\mathcal F} = 0.$$
The solution is a constant. The Lemma follows from 
definition \eqref{lqde}.
\end{proof}

By Lemma \ref{cllw}, we may express the  matrix elements
of ${\mathsf S}$ as ratios of Donaldson-Thomas rubber series:

\begin{equation}\label{lqdet}
\Big\langle \mu \ \Big|
{\mathsf S} \ \Big| \nu \Big \rangle_{\mathcal F}
 = q^{-d}\frac{\Big \langle \mu \ \Big| \frac{1}{1-\psi_\infty}
\ \Big| \nu \Big \rangle^\sim}
{\Big \langle \emptyset \ \Big| \frac{1}{1-\psi_\infty}
\ \Big| \emptyset \Big \rangle^\sim}
.
\end{equation}

\subsection{Computation of the 1-legged vertex}

Let ${\bW'}(\mu,\emptyset,\emptyset)$ denote the
reduced 1-legged vertex,
$${\bW'}(\mu,\emptyset,\emptyset) =
\frac{{\bW}(\mu,\emptyset,\emptyset)}
{{\bW}(\mu,\emptyset,\emptyset)}.$$
The localization formula \eqref{localiz},
\begin{multline*}
\Big \langle\   \Big| \ \lambda([L]) \Big \rangle'_{(-1,0)}=\\ 
 \sum_{|\mu|=d} 
{\bW'}(\mu,\emptyset,\emptyset)|_{s,t_1-s,t_2}
  \cdot 
q^{n(\mu)}  
\bE^{(-1,0)}(\mu) 
\cdot 
\frac{\Big \langle \left[\cI_\mu \right]    \Big| \frac{1}{-s-\psi_\infty}
\ \Big| \lambda \Big \rangle^\sim}
{\Big \langle \emptyset    \Big| \frac{1}{-s-\psi_\infty}
\ \Big| \emptyset \Big \rangle^\sim}
\ t_2^{\ell(\lambda)}, 
\end{multline*}
has a fixed evaluation given by Proposition \ref{dddr} independent of the
equivariant parameters.
If the matrix 
\begin{equation}\label{fth6}
\bE^{(-1,0)}(\mu) 
\cdot 
\frac{\Big \langle \left[\cI_\mu \right]    \Big| \frac{1}{-s-\psi_\infty}
\ \Big| \lambda \Big \rangle^\sim}
{\Big \langle \emptyset    \Big| \frac{1}{-s-\psi_\infty}
\ \Big| \emptyset \Big \rangle^\sim}
\end{equation}
can be calculated, the above localization formula
may be viewed as a square system of 
linear equations for the unknown vector 
${\bW'}(\mu,\emptyset,\emptyset)$,
where $\mu$ ranges over all partitions of a fixed size. 

The matrix \eqref{fth6} is a product of
two factors. The first is an invertible 
diagonal matrix of edge weights explicitly determined by 
\eqref{transf_E}. The second is the operator $\bS$ written 
on the left in the fixed point basis.
The operator $\bS$
is completely determined by the linear differential equation \eqref{xxqq}. 
The inverse of $\bS$ satisfies
$$
q\frac{d}{dq} \bS^{-1} = \MM(0)\, \bS^{-1} - \bS^{-1} \, \MM \,.
$$
In particular, $\bS$ is invertible for $q$ not equal to a root of 
unity.

An identical argument can be used in the local Gromov-Witten theory 
of curves to determine 1-partition Hodge integrals
from the results of \cite{bp}. For Calabi-Yau Hodge integrals, 
the results specialize to the Gopakumar-Mari\~no-Vafa formula
proven in \cite{llz,opgmv}.

Alternatively, both the 1-legged vertex and the 1-partition 
Hodge integrals can be recovered from an parallel localization
formula for the level $(0,0)$ cap.
We leave the details to the reader.


\vspace{+10 pt}
\noindent
Department of Mathematics \\
Princeton University \\
Princeton, NJ 08544, USA\\
okounkov@math.princeton.edu \\

\vspace{+10 pt}
\noindent
Department of Mathematics\\
Princeton University\\
Princeton, NJ 08544, USA\\
rahulp@math.princeton.edu

\end{document}